\documentclass[hidelinks,onefignum,onetabnum]{siamart220329}

\usepackage{amsmath}
\usepackage{amssymb}
\usepackage{bm}
\usepackage{bbm}
\usepackage{mathtools}
\usepackage{adjustbox}
\usepackage{float}
\usepackage{calc}
\usepackage{array}
\usepackage{subcaption}
\usepackage[lighttt]{lmodern}
\usepackage[inline]{enumitem}
\usepackage[font=small]{caption}
\usepackage{algorithm}
\usepackage{algorithmicx}
\usepackage{algpseudocode}
\usepackage{physics}
\usepackage{stmaryrd}
\usepackage{stackengine}
\usepackage{booktabs}
\usepackage{multirow}
\usepackage{threeparttable}
\usepackage{siunitx}
\usepackage{todonotes}
\usepackage{placeins}
\newcommand{\trasSpeedupA}{18.2}
\newcommand{\trasSpeedupB}{39.1}

\newcommand{\trasSolveRate}{98.9\%}
\newcommand{\gdSolveRate}{3.3\%}
\newcommand{\bundleSolveRate}{49.3\%}

\makeatletter
\newcommand\notsotiny{\@setfontsize\notsotiny\@vipt\@viipt}
\makeatother

\newcommand{\mypara}[1]{\noindent\textbf{#1}}

\newcommand{\bigO}[1]{\mathcal{O}\qty(#1)}
\newcommand{\bigOmega}[1]{\operatorname{\Omega}\qty(#1)}
\DeclareMathOperator{\setCl}{cl}
\DeclareMathOperator{\setConv}{conv}
\newcommand{\defeq}{\equiv}
\newcommand{\condSet}[2]{\left\{#1 \;\middle|\; #2\right\}}
\DeclareMathOperator*{\argmax}{argmax}
\DeclareMathOperator*{\argmin}{argmin}
\DeclareMathOperator{\sign}{sign}

\newcommand{\intInt}[2]{\left\llbracket #1,\,#2 \right\rrbracket}   

\newcommand{\V}[1]{\bm{#1}}
\newcommand{\T}[2][]{#2^{#1\intercal}}
\newcommand{\normal}[2]{\mathcal{N}\qty(#1,\,#2)}
\newcommand{\bernoulli}[2]{\mathcal{B}\qty(#1;\,#2)}
\newcommand{\simplexUniform}[1]{\mathcal{U}\qty(\Delta_{#1})}
\newcommand{\prob}[1]{\mathbb{P}\qty[#1]}

\renewcommand{\va}{\V{a}}
\newcommand{\vA}{\V{A}}

\renewcommand{\vb}{\V{b}}
\newcommand{\vB}{\V{B}}
\newcommand{\vc}{\V{c}}
\newcommand{\vC}{\V{C}}
\newcommand{\vd}{\V{d}}

\newcommand{\vg}{\V{g}}
\newcommand{\vG}{\V{G}}
\newcommand{\vh}{\V{h}}
\newcommand{\vI}{\V{I}}
\newcommand{\vl}{\V{l}}
\newcommand{\vM}{\V{M}}

\newcommand{\vp}{\V{p}}

\renewcommand{\vu}{\V{u}}

\newcommand{\vv}{\V{v}}

\newcommand{\vx}{\V{x}}
\newcommand{\vy}{\V{y}}
\newcommand{\vz}{\V{z}}
\newcommand{\vX}{\V{X}}

\renewcommand{\real}{\mathbb{R}}

\newcommand{\nat}{\mathbb{N}}

\newcommand{\posReal}{\mathbb{R}_+}
\newcommand{\nonNegReal}{\mathbb{R}_{\geq 0}}
\newcommand{\integer}{\mathbb{Z}}
\newcommand{\posInt}{\integer_+}

\newcommand{\weaklySmooth}{sparsely nonsmooth}
\newcommand{\betaWeaklySmooth}[1][\beta]{$#1$-sparsely nonsmooth}
\newcommand{\weakSmoothness}{sparse nonsmoothness}
\newcommand{\WeakSmoothness}{Sparse nonsmoothness}

\newcommand{\trasFull}{Trust Region Adversarial Functional Subdifferential}
\newcommand{\tras}{TRAFS}
\algnewcommand\algorithmicinput{\textbf{Input:}}
\algnewcommand\Input{\item[\algorithmicinput]}
\makeatletter
\NewDocumentCommand{\LeftComment}{s m}{%
  \Statex \hspace*{\ALG@thistlm}\(\triangleright\) #2}
\makeatother
\newcommand{\lineRef}[1]{\hyperref[#1]{Line~\ref*{#1}}}
\newcommand{\epsDk}[1][k]{\epsilon_{#1}}    
\newcommand{\subdGenSym}[2]{\partial_{#1}^{\operatorname{#2}}}
\newcommand{\clkSubdSym}[2][\epsilon]{\subdGenSym{#1}{C}#2}
\newcommand{\clkSubd}[3][\epsilon]{\clkSubdSym[#1]{#2}\qty(#3)}
\newcommand{\golSubdSym}[2][\epsilon]{\subdGenSym{#1}{G}#2}
\newcommand{\golSubd}[3][\epsilon]{\golSubdSym[#1]{#2}\qty(#3)}
\newcommand{\relSubdSym}[2][\epsilon]{\subdGenSym{#1}{r}#2}
\newcommand{\relSubd}[3][\epsilon]{\relSubdSym[#1]{#2}\qty(#3;\,S)}

\newcommand{\funSubdGSym}[2][\epsilon]{\partial_{#1}#2}
\newcommand{\funSubdDSym}[2][\epsilon]{\mathring{D}_{#1}#2}
\newcommand{\funSubdG}[3][\epsilon]{\funSubdGSym[#1]{#2}\qty(#3)}
\newcommand{\funSubdGE}[2]{\funSubdG[\epsDk]{#1}{#2}}
\newcommand{\funSubdD}[3][\epsilon]{\funSubdDSym[#1]{#2}\qty(#3)}
\newcommand{\funSubdP}[2][\epsilon]%
{\qty(\funSubdDSym[#1]{#2},\,\funSubdGSym[#1]{#2})}

\newcommand{\fmeth}{f_{\text{meth}}}
\newcommand{\fopt}{f_{\text{opt}}}

\newcommand{\SAtwo}{$\text{SA}_2$}

\newcommand{\mytitle}{\tras{}: A Nonsmooth Convex Optimization Algorithm with
$\bigO{\frac{1}{\epsilon}}$ Iteration Complexity}
\newcommand{\shorttitle}{\tras{}: A Nonsmooth Convex Optimization Algorithm}
\newcommand{\shortauthors}{K. Jia, M. Rinard}

\headers{\shorttitle}{\shortauthors}

\title{\mytitle}

\author{Kai Jia\thanks{MIT CSAIL, Cambridge, MA (\email{jiakai@mit.edu})}
\and Martin Rinard\thanks{MIT CSAIL, Cambridge, MA (\email{rinard@csail.mit.edu})}
}

\newsiamremark{remark}{Remark}

\usepackage[numbers]{natbib}

\ifpdf
\hypersetup{
  pdftitle={\shorttitle},
  pdfauthor={\shortauthors}
}
\fi

\begin{document}

\maketitle
\begin{abstract}
    We present the \trasFull{} (\tras) algorithm for constrained optimization of
    nonsmooth convex Lipschitz functions. Unlike previous methods that assume a
    subgradient oracle model, we work with the \emph{functional subdifferential}
    defined as a set of subgradients that simultaneously captures sufficient
    local information for effective minimization while being easy to compute for
    a wide range of functions. In each iteration, \tras{} finds the best step
    vector in an $\ell_2$-bounded trust region by considering the worst bound
    given by the functional subdifferential. \tras{} finds an approximate
    solution with an absolute error up to~$\epsilon$ in $\bigO{ \epsilon^{-1}}$
    or $\bigO{\epsilon^{-0.5}}$ iterations depending on whether the objective
    function is strongly convex, compared to the previously best-known bounds of
    $\bigO{\epsilon^{-2}}$ and $\bigO{\epsilon^{-1}}$ in these settings. \tras{}
    makes faster progress if the functional subdifferential satisfies a locally
    quadratic property; as a corollary, \tras{} achieves linear convergence
    (i.e., $\bigO{\log \epsilon^{-1}}$) for strongly convex smooth functions. In
    the numerical experiments, \tras{} is on average \trasSpeedupB{}x faster and
    solves twice as many problems compared to the second-best method.
\end{abstract}



\section{Introduction}
We consider the problem $\min_{\vx \in C} f(\vx)$ where $f: \real^n \mapsto
\real$ is a convex function, and the constraint set $C \subset \real^n$ is
compact and convex. The function $f(\cdot)$ is Lipschitz continuous over $C$ but
may be nondifferentiable on certain parts of~$C$.

There is a long history of research on such problems~\citep{
lemarechal1978nonsmooth, hiriart1996convex, bagirov2014introduction,
bagirov2020numerical}. Most existing methods assume an oracle that delivers an
arbitrary subgradient of $f(\cdot)$ at any given point. Such methods iteratively
update the solution based on the current and previous subgradients. However, an
arbitrary subgradient does not contain sufficient information to characterize
the local behavior of the objective function; following the subgradient
direction may even increase the objective value. Consequently, such methods need
at least $\bigOmega{\epsilon^{-2}}$ iterations to find an $\epsilon$-approximate
solution in the worst case if the oracle is queried once per iteration \citep{
nemirovskij1983problem, iouditski2014primal}.

Instead of working with an arbitrary subgradient, we propose the
\emph{functional subdifferential} (\cref{def:tras:func-subd}), which is a set of
subgradients that incorporates sufficient local information for effective
optimization while being easy to compute for a large class of functions of
interest (\cref{rmk:tras:func-subd:easy-to-compute}). The \emph{\trasFull}
(\tras) algorithm (\cref{algo:tras}) utilizes the functional subdifferential for
nonsmooth optimization. In each iteration, \tras{} chooses a step vector by
solving a minimax problem involving an $\ell_2$-bounded trust region and the
functional subdifferential. Since the minimax problem adversarially considers
the worst bound of the change of the objective value, \tras{} guarantees that
the objective value never increases at any iteration.

For Lipschitz functions, \tras{} achieves an iteration complexity of
$\bigO{\epsilon^{-1}}$ (\cref{thm:tras:rate}). If the function is also strongly
convex, \tras{} converges in $\bigO{\epsilon^{-0.5}}$ iterations
(\cref{thm:tras:sc:rate}). These results improve on the previously best-known
iteration complexities of $\bigO{ \epsilon^{-2}}$ for general Lipschitz
functions and $\bigO{\epsilon^{-1}}$ for strongly convex Lipschitz
functions~\citep{bubeck2015convex, diaz2023optimal}. Because the functional
subdifferential contains multiple subgradients, we escape the corresponding
lower bounds of $\bigOmega{ \epsilon^{-2}}$ and $\bigOmega{\epsilon^{-1}}$ for
optimization algorithms that use a single subgradient per iteration~\citep{
nemirovskij1983problem, iouditski2014primal}.

When the functional subdifferential satisfies the locally quadratic property
(\cref{def:tras:quad-func-subd}), which intuitively means that it incorporates
information from a quadratically larger neighborhood, \tras{} makes faster
iterate progress (\cref{thm:tras:scqf:rate}). As a corollary, \tras{} achieves
linear convergence (i.e., $\bigO{\log \epsilon^{-1}}$) for strongly convex
smooth functions (\cref{thm:tras:sc-smooth:rate}), which recovers the lower
bound of first-order methods for this class of
functions~\citep{bubeck2015convex}.

The above complexity results of \tras{} assume the ability to solve a minimax
problem involving the functional subdifferential in each iteration. Although the
functional subdifferential can be defined for all Lipschitz convex functions
(\cref{thm:tras:func-subd:from-clarke,thm:tras:func-subd:from-goldstein}), the
minimax problem with universally-defined functional subdifferentials is
typically intractable. Fortunately, a function usually has more than one
functional subdifferential, with some easier to work with than others. We
present compositional rules to compute the functional subdifferential for
functions that use common operators (such as sum, max reduction, linear
transform, etc.; see \cref{sec:func-subd-rule}). Our set of rules enables
efficient implementation of \tras{} for a wide range of functions, including
nonsmooth convex benchmark problems used in previous work~\citep{
haarala2004new, bagirov2014introduction, nesterov2015quasi, karmitsa2020limited}
and various hard-case functions constructed in the proofs of the aforementioned
lower bounds of iteration complexities \citep{ bubeck2015convex,
iouditski2014primal, agarwal2009information}. When an efficient functional
subdifferential is not available, we use the \emph{almost-functional
subdifferential} with weaker properties (\cref{def:tras:almost-func-subd}) at
the cost of higher iteration complexities of $\bigO{\epsilon^{-2}}$ and
$\bigO{\epsilon^{-1}}$ (\cref{thm:tras:almost-func-subd-rate}).

We present an adaptive \tras{} implementation using heuristics inspired by our
convergence analysis, which requires no knowledge of the convexity or smoothness
parameters. In numerical experiments consisting of benchmark problems used in
previous work \citep{ haarala2004new, bagirov2014introduction,
nesterov2015quasi, karmitsa2020limited} and new problems proposed in this paper,
\tras{} is on average \trasSpeedupB{} times faster and solves twice as many
problems compared to the second-best method
(\cref{sec:tras:experiments:results}).

\section{Related work}
\label{sec:tras:related}
There is a large body of work on smooth and nonsmooth convex optimization
\citep{boyd2004convex, bubeck2015convex, bagirov2020numerical}. We review the
most relevant work for unstructured nonsmooth convex optimization.

We review two classic methods that use the subgradient oracle. The projected
subgradient descent computes $\vx_{k+1} = \Pi_C\qty(\vx_k - \eta_k \vg_k)$,
where $\eta_k \in \posReal$ is the step size, $\vg_k \in \partial f(\vx_k)$ is
an arbitrary subgradient, and $\Pi_{C}(\vx) \defeq \argmin_{\vy \in C} \norm{\vy
- \vx}_2$ is the projection operator onto $C$. Setting $\eta_k =
\frac{R}{L\sqrt{k}}$ yields an ergodic convergence rate of $\bigO{
\epsilon^{-2}}$ \citep{bubeck2015convex}. The bundle method
\citep{lemarechal1978nonsmooth, mifflin1982modification, hiriart1996convex}
approximates the objective at the $k$-th iteration by a piecewise affine
function $\tilde{f}_k(\vx) \defeq \max_{j \in J_k} \qty( f(\vx_j) +
\T{\vg_j}(\vx - \vx_j))$ where $J_k \subset \intInt{1}{k}$. Different variants
may employ different strategies of defining $J_k$ or computing $\vx_{k+1}$ from
$\tilde{f}_k(\cdot)$~\citep{lemarechal1995new, kiwiel2006methods,
makela2016proximal}. It was recently shown that the proximal bundle method has
iteration complexities of $\bigO{\epsilon^{-2}}$ or $\bigO{\epsilon^{-1}}$ for
Lipschitz and strongly convex Lipschitz functions, respectively~\citep{
diaz2023optimal}. Compared to these two methods, \tras{} provides both better
iteration complexity guarantees and substantially better performance on our
benchmark problems (\cref{sec:tras:experiments:results}).

Any method that queries a subgradient oracle once per iteration (thus including
the projected subgradient descent and the bundle method) needs at least
$\bigOmega{\epsilon^{-2}}$ or $\bigOmega{\epsilon^{-1}}$ iterations in the worst
case, depending on whether the function is strongly
convex~\citep{nemirovskij1983problem, iouditski2014primal}. The same lower
bounds also hold in the stochastic setting where the oracle delivers a
subgradient with some zero-mean, bounded-variance additive noise, which can be
proven from an information-theoretic perspective \citep{
agarwal2009information}. By contrast, \tras{} uses the functional
subdifferential that includes multiple subgradients and incorporates sufficient
local information to achieve a lower iteration complexity.

There are attempts to utilize the $\epsilon$-subdifferential (see
\cref{rmk:tras:relaxed-subd} for a definition) for optimizing nonsmooth
functions with special structures \citep{bertsekas1971steepest,
cullum1975minimization}. Although rules exist to characterize the
$\epsilon$-subdifferentials mathematically~\citep{ hiriart1996convex}, computing
the $\epsilon$-subdifferential for unstructured functions is typically
intractable \citep{ burke2002approximating}. Gradient sampling methods
\citep{burke2002approximating, burke2005robust, burke2020gradient} overcome this
intractability to optimize unstructured functions by using gradients computed at
$m$ points uniformly sampled in $B_{\epsilon_k}[\vx_k]$ in each iteration to
approximate the $\epsilon$-subdifferential. However, gradient sampling methods
are computationally expensive since most variants require $m \geq n+1$,
where~$n$ is the dimension of the problem. Even the variant converging with a
constant number of samples still recommends $m = n/10$ for practical performance
\citep{ curtis2013adaptive}. Previous evaluations typically work with
small-scale problems with $n \leq 120$ \citep{ burke2005robust,
curtis2013adaptive, helou2017local}. We are unaware of any global iteration
complexity results for gradient sampling methods other than a local linear
convergence result under a special class of functions~\citep{ helou2017local}.
By contrast, the functional subdifferential deterministically characterizes the
local behavior of the function to enable the guaranteed global convergence rates
of \tras{}.

\section{Notation}
We use bold symbols to denote vectors (lower case) and matrices~(upper case).
For example, $\va$ is a vector, $\va_i$ is a vector indexed by $i$, $a_i$ is the
$i$-th element of $\va$, $\vA_i$ is the $i$-th row of the matrix $\vA$, $\V{1}$
is the all-one vector whose dimension is inferred from context, and $\V{1}_n$ is
the all-one vector of dimension $n$. We use $k\Delta_n \defeq \condSet{\vx \in
\real^n}{\T{\V{1}}\vx = k,\, \vx \geq \V{0}}$ to denote the $n$-dimensional
simplex scaled by $k$ ($k$ is omitted if $k=1$), $B_{r}[\vx] \defeq \condSet{\vy
\in \real^n}{\norm{\vy - \vx}_2 \leq r}$ for the $\ell_2$-ball centered at $\vx$
with radius $r$, $\intInt{a}{b} \defeq [a,\, b] \cap \integer$ for the set of
integers between $a$ and $b$, $\setCl A$ for the closure of a set $A$, and
$\setConv A$ for the convex hull of a set $A$.

\section{The functional subdifferential}
\subsection{Introducing the functional subdifferential}
We first recall the definitions of subgradient and subdifferential. We then
define the \emph{relaxed subdifferential} as a natural idea to include the
subdifferential of local nonsmooth points. We define a class of nonsmooth
functions called the \emph{\weaklySmooth} functions. We show that the relaxed
subdifferential provides a local upper bound and global lower bound for
\weaklySmooth{} functions. Finally, we generalize the definition of the relaxed
subdifferential to the \emph{functional subdifferential} by only requiring the
same upper and lower bounds to be met.

\begin{definition}[{{Subdifferential
    \citep[Definitions 3.1,3.2]{beck2017first}}}]
    Let $f: \real^n \mapsto \real$ be a convex function. A vector $\vg \in
    \real^n$ is a \emph{subgradient} of $f(\cdot)$ at $\vx \in \real^n$ if and
    only if $f(\vy) \geq f(\vx) + \T{\vg}\qty(\vy - \vx)$ holds for all $\vy \in
    \real^n$. The set of all subgradients of $f(\cdot)$ at $\vx$ is called the
    \emph{subdifferential} of $f(\cdot)$ at $\vx$ and is denoted by $\partial
    f(\vx)$:
    \begin{align}
        \partial f(\vx) \defeq \condSet{\vg \in \real^n}{
            \forall \vy \in \real^n:\: f(\vy) \geq f(\vx) +
            \T{\vg}\qty(\vy - \vx)}
    \end{align}
\end{definition}

\begin{proposition}
    \label{thm:subgrad-prop}
    The subdifferential of a convex function $f: \real^n \mapsto \real$
    satisfies a few properties~\citep[Chapter 3]{beck2017first}:
    \begin{itemize}
        \item For all $\vx \in \real^n$, $\partial f(\vx)$ is nonempty, convex,
            closed, and bounded.
        \item If $f(\cdot)$ is differentiable at $\vx$, then $\partial f(\vx) =
            \{\nabla f(\vx)\}$.
        \item For a convex set $C \subset \real^n$, $\vx^* \in \argmin_{\vx\in
            C}f(\vx)$ if and only if there exists $\vg \in \partial f(\vx^*)$
            such that $\forall \vy \in C:\: \T{\vg}\qty(\vy - \vx^*) \geq 0$.
        \item If $f(\cdot)$ is $L$-Lipschitz over an open set $S$, then
            \begin{align}
                \forall \vx \in S:\:
                \forall \vg \in \partial f(\vx):\:
                \norm{\vg}_2 \leq L
                \label{eqn:tras:relaxed-subg-norm-bound}
            \end{align}
    \end{itemize}
\end{proposition}

\begin{definition}[Relaxed subdifferential]
    \label{def:tras:relaxed-subd}
    For a convex function $f: \real^n \mapsto \real$, a set $S \subset \real^n$,
    and $\epsilon \in \posReal$, its $\epsilon$-relaxed subdifferential at $\vx
    \in \real^n$ constrained to $S$ is defined as
    \begin{align*}
        \relSubd{f}{\vx} &\defeq \setCl \setConv
            \condSet{\vg}{\vg \in \partial f(\vy) \text{ for }
                 \vy = \vx \text{ or }
             \qty(\vy \in S \cap B_\epsilon[\vx] \text{ and }
             \abs{\partial f(\vy)} > 1)}
    \end{align*}
\end{definition}
\begin{remark}
    \label{rmk:tras:relaxed-subd}
    The relaxed subdifferential differs from the $\epsilon$-subdifferential
    (defined as $\clkSubd{f}{\vx} \defeq \condSet{\vg \in \real^n}{\forall \vy
    \in \real^n:\: f(\vy) \geq f(\vx) + \T{\vg}(\vy - \vx) - \epsilon}$, cf.
    \citep[Definition 2.17]{ bagirov2014introduction}) and the Goldstein
    $\epsilon$-subdifferential ($\golSubd{f}{\vx} \defeq \setCl \setConv
    \bigcup_{\vy \in B_{\epsilon}[\vx]} \partial f(\vy)$, cf. \citep[Definition
    3.3]{bagirov2014introduction}) in that the relaxed subdifferential only
    considers the local nondifferentiable points can can be easier to compute.
\end{remark}

Now we introduce \weakSmoothness{} and related properties. Informally, a
\weaklySmooth{} function has finitely separated nondifferentiable points along
any direction, and the gradient over differentiable intervals is Lipschitz
continuous.
\begin{definition}[\WeakSmoothness]
    \label{def:weak-smooth}
    Given a function $f: \real^n \mapsto \real$ and a convex set $S \subset
    \real^n$, $f(\cdot)$ is \emph{\betaWeaklySmooth} over $S$ if and only
    if:

    For all $(\vx_0,\,\vd) \in S \times \real^n$ such that $\norm{\vd}_2 = 1$
    and $\vx_0 + k\vd \in S$ for some $k \neq 0$, define $g(\lambda) \defeq
    f(\vx_0 + \lambda \vd)$; define $E \defeq \condSet{\lambda \in \real}{\vx_0
    + \lambda \vd \in S}$ as the interval of valid values of $\lambda$ regarding
    $S$. Then there exists a (possibly empty or infinite) sequence $(c_1,\,
    \ldots,\, c_m)$ in $E$ such that $\qty(\inf_{i\in \intInt{1}{m-1}} c_{i+1} -
    c_i) > 0$, $g(\cdot)$ is not differentiable at $c_i$ for $i \in
    \intInt{1}{m}$, and $g(\cdot)$ is $\beta$-smooth over all intervals $(c_i,\,
    c_{i+1})$ for $i \in \intInt{0}{m}$, where $c_0 \defeq \inf E$ and $c_{m+1}
    \defeq \sup E$, i.e., for all $i \in \intInt{0}{m}$, for all $(\lambda_1,\,
    \lambda_2) \in (c_i,\, c_{i+1})^2$, we have $\abs{g'(\lambda_2) -
    g'(\lambda_1)} \leq \beta \abs{\lambda_2 - \lambda_1}$.
\end{definition}

\begin{proposition}
    If a function $f(\cdot)$ is $\beta$-smooth over a convex set $S$, then it is
    also \betaWeaklySmooth{} over $S$.
\end{proposition}
\begin{proof}
    Given $\vx_0$ and $\vd$ as in \cref{def:weak-smooth}, we have $g'(\lambda) =
    \T{\vd}\nabla f(\vx_0 + \lambda \vd)$. Therefore,
    \begin{align*}
        \abs{g'(\lambda_2)-g'(\lambda_1)}
            &\leq \norm{\vd}_2 \norm{\nabla f(\vx_0+\lambda_2\vd) -
                \nabla f(\vx_0+\lambda_1\vd)}_2  \\
            &\leq \beta \norm{\qty(\vx_0+\lambda_2\vd) -
                \qty(\vx_0+\lambda_1\vd)}_2
            = \beta\abs{\lambda_2 - \lambda_1}
    \end{align*}
\end{proof}

The relaxed subdifferential is related to the \weakSmoothness{} as the
following:
\begin{lemma}
    \label{thm:tras:relaxed-subg-bounds}
    Given $\epsilon \in \posReal$, a convex function $f: \real^n \mapsto \real$,
    and an open convex set $S \subset \real^n$, if $f(\cdot)$ is $L$-Lipschitz
    over $S$, then for any $(\vx,\,\vy) \in S^2$,
    \begin{align}
        f(\vy) \geq \qty(
            f(\vx) + \max_{\vg \in \relSubd{f}{\vx}} \T{\vg}(\vy-\vx))
            - 2\epsilon L
        \label{eqn:tras:relaxed-subg-lb}
    \end{align}
    Moreover, if $f(\cdot)$ is also \betaWeaklySmooth{} over $S$ and
    $\norm{\vy - \vx}_2 \leq \epsilon$, then
    \begin{align}
        f(\vy) \leq \qty(
            f(\vx) + \max_{\vg \in \relSubd{f}{\vx}} \T{\vg}(\vy-\vx))
            + \frac{\beta}{2} \norm{\vy-\vx}_2^2
        \label{eqn:tras:relaxed-subg-ub}
    \end{align}
\end{lemma}
\begin{proof}
    A proof of \cref{eqn:tras:relaxed-subg-lb} can found at \citet[Theorem
    3.12]{bagirov2014introduction}.

    Before proving \cref{eqn:tras:relaxed-subg-ub}, let's recall the definition
    of directional derivative: for a convex function $h: \real^m \mapsto \real$
    and a vector $\vd \in \real^m$, the directional derivative of $h(\cdot)$ at
    $\vx \in \real^m$ in the direction $\vd$ is defined as $h'(\vx;\vd) \defeq
    \lim_{\alpha \to 0^+} \frac{h(\vx + \alpha \vd) - h(\vx)}{\alpha}$. It is
    related to the subdifferential through $\max_{\vg \in \partial h(\vx)}
    \T{\vg}\vd = h'(\vx;\vd)$ \citep[Theorem 3.26]{beck2017first}. When the
    dimension $m=1$, we denote $h'_+(x) \defeq h'(x;\,1)$ and $h'_-(x) \defeq
    h'(x;\,-1)$. If $m=1$ and $h(\cdot)$ is differentiable over the open
    interval $(x,\, y)$, then $h'_+(x) = \lim_{\alpha \to x^+} h'(\alpha)$ and
    $h'_-(y) = \lim_{\alpha \to y^-} h'(\alpha)$ \citep[Theorem 24.1]{
    rockafellar1970convex}.

    Now let's assume $f(\cdot)$ is \betaWeaklySmooth{} over $S$ and $\norm{\vy
    - \vx}_2 \leq \epsilon$. Let $r \defeq \norm{\vy - \vx}_2$ and $\vd \defeq
    (\vy - \vx) / r$. Define $g(\lambda) \defeq f(\vx + \lambda \vd)$. Due to
    the definition of \weakSmoothness{}, there is a finite sequence $(
    c_0\defeq 0,\, c_1,\, \ldots,\, c_m,\, c_{m+1}=r)$ in $[0,\,r]$ such that
    $g(\cdot)$ is $\beta$-smooth over $(c_i,\, c_{i+1})$ for $i \in \intInt{
    0}{m}$ and $g(\cdot)$ is not differentiable at $c_i$ for $i \in
    \intInt{1}{m}$. Thus we have $\partial f(\vx + c_i \vd) \subset
    \relSubd{f}{\vx}$ for~$i \in \intInt{0}{m}$.

    For each $i \in \intInt{0}{m}$, we define a $\beta$-smooth function
    $g_i(\lambda)$ over $[c_i,\, c_{i+1}]$ such that $g_i'(\lambda) =
    g'(\lambda)$ for $c_i < \lambda < c_{i+1}$, $g_i'(c_i) = g'_+(c_i)$,
    $g_i'(c_{i+1}) = g'_-(c_{i+1})$, and $g_i(c_i) = g(c_i)$. It can be easily
    shown $g_i(\lambda) = g(\lambda)$ for $c_i \leq \lambda \leq c_{i+1}$.

    Let $M \defeq \max_{\vg \in \relSubd{f}{\vx}} \T{\vg}(\vy - \vx)$. Let $s_i
    \defeq c_{i+1} - c_i$. For $i \in \intInt{0}{m}$, we have $g_i'(c_i) =
    g'_+(c_i) = f'(\vx + c_i \vd; \vd) = \max_{\vg \in \partial f(\vx + c_i\vd)}
    \T{\vg}\vd \leq \max_{\vg \in \relSubd{f}{\vx}} \T{\vg}\vd = \frac{M}{r}$.
    The $\beta$-smoothness of $g_i(\cdot)$ implies that $g_i(c_{i+1}) - g_i(c_i)
    \leq g_i'(c_i) (c_{i+1} - c_i) + \frac{\beta}{2} (c_{i+1} - c_i)^2 \leq
    \frac{M}{r} s_i + \frac{\beta}{2} s_i^2$. Therefore,
    \begin{align*}
        f(\vy) - f(\vx) &= g(r) - g(0) =
            \sum_{i=0}^m \qty(g_i(c_{i+1}) - g_i(c_i))
            \leq \sum_{i=0}^m \qty(\frac{M}{r}s_i + \frac{\beta}{2}s_i^2) \\
            & \leq \frac{M}{r} \qty(\sum_{i=0}^M s_i) +
                \frac{\beta}{2}\qty(\sum_{i=0}^m s_i )^2
            = \max_{\vg \in \relSubd{f}{\vx}} \T{\vg}(\vy - \vx) +
                \frac{\beta}{2} \norm{\vy-\vx}_2^2
    \end{align*}
\end{proof}

\Cref{eqn:tras:relaxed-subg-lb,eqn:tras:relaxed-subg-ub} are the core properties
that would enable the convergence rates of the proposed \tras{} algorithm. One
crucial observation is that the relaxed subdifferential is not the only mapping
that satisfies these properties. Therefore, we propose the following
\emph{functional subdifferential} that captures any mapping satisfying
\cref{eqn:tras:relaxed-subg-lb,eqn:tras:relaxed-subg-ub}.

\begin{definition}[Functional subdifferential]
    \label{def:tras:func-subd}
    Given a convex function $f: \real^n \mapsto \real$ and an open convex set
    $S \subset \real^n$, a pair $\funSubdP{f}$ is called a \emph{functional
    subdifferential} of $f(\cdot)$ over $S$ where $\epsilon \in \nonNegReal$
    (called the \emph{slack}), $\funSubdGSym{f}: S \mapsto 2^{\real^n}$, and
    $\funSubdDSym{f}: S \mapsto \nonNegReal \cup \qty{ +\infty}$, if and only if
    the following properties hold:
    \begin{itemize}
        \item For any $\vx \in S$, $\funSubdG{f}{\vx}$ is a nonempty, convex,
            closed, and bounded set.
        \item There exists $L \in \posReal \cup \qty{0^+}$ such that
            $\funSubdD{f}{\vx} \geq \frac{\epsilon}{2L}$ for $\vx \in S$ and
            $\epsilon \in \posReal$.
        \item For any $\vx \in S$, $\epsilon \in \nonNegReal$, and $\vy \in S$,
            it holds that
            \begin{align}
                f(\vy) \geq
                    \qty(f(\vx) +
                        \max_{\vg \in \funSubdG{f}{\vx}} \T{\vg}(\vy - \vx)) -
                \epsilon
                \label{eqn:tras:func-subd:lb}
            \end{align}
        \item There exists $\beta \in \nonNegReal$ such that for any $\vx \in
            S$, $\epsilon \in \nonNegReal$, and $\vy \in S$ such that $\norm{\vy
            - \vx}_2 \leq \funSubdD{f}{\vx}$, it holds that
            \begin{align}
                f(\vy) \leq
                    \qty(f(\vx) +
                        \max_{\vg \in \funSubdG{f}{\vx}} \T{\vg}(\vy - \vx)) +
                    \frac{\beta}{2} \norm{\vy - \vx}_2^2
                \label{eqn:tras:func-subd:ub}
            \end{align}
    \end{itemize}
    The constants $L$ and $\beta$ are called the \emph{associated constants} of
    the functional subdifferential $\funSubdP{f}$. Of note, for $L' \geq L$ and
    $\beta' \geq \beta$, $(L',\, \beta')$ is also a pair of associated constants
    of $\funSubdP{f}$.
\end{definition}
\begin{remark}
    We have included $\funSubdD{f}{\vx}$ in the functional subdifferential
    definition to simplify the presentation and analysis; alternatively, we
    could define $\funSubdD{f}{\vx}$ as the maximum distance between $\vx$ and
    $\vy$ such that \cref{eqn:tras:func-subd:ub} holds. The \tras{} algorithm
    only needs the solution of a minimax problem involving $\funSubdG{f}{\vx}$.
\end{remark}

\subsection{Rules for computing the functional subdifferential}
\label{sec:func-subd-rule}
This subsection presents compositional rules to compute the functional
subdifferential. We assume all functions are convex in $\real^n$. For a function
$f(\cdot)$, we denote its functional subdifferential as $\funSubdP{f}$ and the
associated constants as $\qty(L_f,\, \beta_f)$. We assume all functional
subdifferentials are defined over an open convex set $S \subset \real^n$ unless
otherwise specified.

We start with four terminal cases for constructing the functional
subdifferential.

\begin{proposition}[Functional subdifferential for smooth functions]
    \label{thm:tras:func-subd:from-smooth}
    If $f(\cdot)$ is $\beta$-smooth over $S$, then $\funSubdP{f}$ is a
    functional subdifferential of~$f(\cdot)$ with associated constants $(0^+,\,
    \beta)$ where
    \begin{align*}
        \funSubdD{f}{\vx} \defeq +\infty,\quad
        \funSubdG{f}{\vx} \defeq \qty{\nabla f(\vx)}
    \end{align*}
\end{proposition}
\begin{proof}
    It is straightforward to verify the properties in \cref{def:tras:func-subd}.
\end{proof}

\begin{proposition}[Functional subdifferential for \weaklySmooth{} functions]
    \label{thm:tras:func-subd:from-relaxed}
    Assume $f(\cdot)$ is $L$-Lipschitz and \betaWeaklySmooth{} over $S$.
    Let $\relSubd{f}{\vx}$ be the relaxed subdifferential of $f(\cdot)$ as
    defined in \cref{def:tras:relaxed-subd}. Then $\funSubdP{f}$ is a functional
    subdifferential of $f(\cdot)$ with associated constants $(L,\, \beta)$ where
    \begin{align*}
        \funSubdD{f}{\vx} \defeq \frac{\epsilon}{2L},\quad
        \funSubdG{f}{\vx} \defeq \relSubd[\frac{\epsilon}{2L}]{f}{\vx}
    \end{align*}
\end{proposition}
\begin{proof}
    Use \cref{thm:tras:relaxed-subg-bounds} to verify the properties in
    \cref{def:tras:func-subd}.
\end{proof}

\begin{proposition}[Functional subdifferential from the
    $\epsilon$-subdifferential]
    \label{thm:tras:func-subd:from-clarke}
    If $f(\cdot)$ is $L$-Lipschitz over $S$, then $\funSubdP{f}$ is a functional
    subdifferential of~$f(\cdot)$ with associated constants $(L,\, 0)$ where
    \begin{align*}
        \funSubdD{f}{\vx} &\defeq \frac{\epsilon}{2L},\quad
        \funSubdG{f}{\vx} \defeq \condSet{\vg \in \real^n}{\forall \vy
        \in \real^n:\: f(\vy) \geq f(\vx) + \T{\vg}(\vy - \vx) - \epsilon}
    \end{align*}
\end{proposition}
\begin{proof}
    A proof of the first property in \cref{def:tras:func-subd} can found at
    \citet[Theorem 2.32]{bagirov2014introduction}. \Cref{eqn:tras:func-subd:lb}
    holds by the definition of $\funSubdG{f}{\vx}$. For $(\vx,\,\vy,\,\epsilon)
    \in S^2\times \nonNegReal$ such that $\norm{\vy - \vx}_2 \leq
    \funSubdD{f}{\vx}$, we have $\partial f(\vy) \subset \funSubdG{f}{\vx}$
    \citep[Theorem 2.33]{bagirov2014introduction}. For any $\vg_y \in \partial
    f(\vy)$, we have $f(\vx) \geq f(\vy) + \T{\vg_y}(\vx - \vy)$, which implies
    $\max_{\vg \in \funSubdG{f}{\vx}} \T{\vg}(\vy - \vx) \geq \T{\vg_y}(\vy -
    \vx) \geq f(\vy) - f(\vx)$ and thus proves~\cref{eqn:tras:func-subd:ub}.
\end{proof}

\begin{proposition}[Functional subdifferential from the Goldstein
    subdifferential]
    \label{thm:tras:func-subd:from-goldstein}
    If $f(\cdot)$ is $L$-Lipschitz over $S$, then $\funSubdP{f}$ is a functional
    subdifferential of~$f(\cdot)$ with associated constants $(L,\, 0)$ where
    \begin{align*}
        \funSubdD{f}{\vx} &\defeq \frac{\epsilon}{2L},\quad
        \funSubdG{f}{\vx} \defeq \golSubd[\frac{\epsilon}{2L}]{f}{\vx},\quad
        \golSubd{f}{\vx} \defeq \setCl \setConv
            \bigcup_{\vy \in B_{\epsilon}[\vx]} \partial f(\vy)
    \end{align*}
\end{proposition}
\begin{proof}
    A proof of \cref{eqn:tras:func-subd:lb} can found at \citet[Theorem
    3.12]{bagirov2014introduction}. \Cref{eqn:tras:func-subd:ub} can be
    proven similarly to \cref{thm:tras:func-subd:from-clarke} since $ f(\vy)
    \subset \funSubdG{f}{\vx}$ by definition.
\end{proof}

Next we present a few compositional rules for computing the functional
subdifferential through common operations.

\begingroup
\newcommand{\hb}[1][\epsilon]{\bar{h'}_{#1}(\vx)}
\begin{proposition}[Rule of composition]
    \label{thm:tras:func-subd:outer-comp}
    Let $F(\vx) \defeq h(f(\vx))$ where $f: \real^n \mapsto \real$ is a convex
    function and $h: \real \mapsto \real$ is a convex non-decreasing function.
    Let $S \subset \real^n$ be an open convex set. Assume $f(\cdot)$ is
    $L_1$-Lipschitz over $S$ and $\funSubdP{f}$ is a functional subdifferential
    of $f(\cdot)$ over $S$ with associated constants $(L_f,\, \beta_f)$. Let $T
    \defeq \condSet{f(\vx)}{\vx \in S}$. Assume $h(\cdot)$ is $L_h$-Lipschitz
    and \betaWeaklySmooth[\beta_h]{} over $T$ and $\funSubdP{h}$ is the
    functional subdifferential of $h(\cdot)$ over $T$ defined by
    \cref{thm:tras:func-subd:from-relaxed} or
    \cref{thm:tras:func-subd:from-smooth}. Define
    \begin{align}
        \xi(\vx,\,\epsilon) &\defeq \max\condSet{
        \xi \in \qty[0,\, \frac{\epsilon}{\hb}]}{
            \substack{
                \forall \vy \in S:\: \norm{\vy - \vx}_2 \leq \xi \implies \\
                \abs{f(\vy) -  f(\vx)} \leq \funSubdD[\epsilon -
                \hb\xi]{h}{f(\vx)}
        }}
        \label{eqn:tras:func-subd:outer-comp:xi-def}
        \\
        \text{where } \hb &\defeq \sup_{\epsilon'\in[0,\,\epsilon]}
            \max \funSubdG[\epsilon']{h}{f(\vx)} \nonumber
    \end{align}
    Then $F(\cdot)$ is a convex function that has a functional subdifferential
    $\funSubdP{F}$ with associated constants $(L_F,\, \beta_F)$ where
    \begin{align*}
        \funSubdD{F}{\vx} &\defeq \funSubdD[\xi(\vx,\,\epsilon)]{f}{\vx} \\
        \funSubdG{F}{\vx} &\defeq \condSet{
            \alpha \vg}{
            \alpha \in \funSubdG[\gamma(\vx,\,\epsilon)]{h}{f(\vx)},\:
            \vg \in \funSubdG[\xi(\vx,\,\epsilon)]{f}{\vx}} \\
        \gamma(\vx,\,\epsilon) &\defeq \epsilon - \hb\xi(\vx,\,\epsilon) \\
        L_F \leq L_t L_f,\quad
        &L_t \defeq L_h(2L_1 + 1),\quad
        \beta_F \defeq L_h\beta_f + L_1^2\beta_h
    \end{align*}
    Alternatively, one can set $\xi(\vx,\,\epsilon) \defeq \frac{\epsilon}{
    L_t}$, which may yield smaller $\funSubdD{F}{\vx}$.
\end{proposition}
\begin{proof}
    It is straightforward to verify that $F(\cdot)$ is convex and that the first
    property in \cref{def:tras:func-subd} is satisfied. Of note,
    $\funSubdG{h}{\cdot}$ is a closed interval, and we have $0 \leq \hb \leq
    L_h$ due to $h(\cdot)$ being non-decreasing and
    \cref{eqn:tras:relaxed-subg-norm-bound}. We have $\frac{\epsilon}{L_t} \in
    \qty[0,\, \frac{\epsilon}{\hb}]$ since $\frac{\epsilon}{L_t} \leq
    \frac{\epsilon}{L_h} \leq \frac{\epsilon}{\hb}$. We also have $\funSubdD[
    \epsilon - \hb \frac{\epsilon}{L_t}]{h}{f(\vx)} \geq
    \frac{1}{2L_h}\qty(\epsilon - \hb \frac{\epsilon}{L_t}) \geq
    \frac{\epsilon}{2L_h} \qty(1 - \frac{L_h}{L_t}) = \frac{\epsilon
    L_1}{L_h(2L_1+1)} = L_1 \frac{\epsilon}{L_t} \geq \abs{f(\vy)-f(\vx)}$ when
    $\norm{\vy-\vx}_2 \leq \frac{\epsilon}{L_t}$. Combining these two facts
    yields $\xi(\vx,\, \epsilon) \geq \frac{\epsilon}{L_t}$ and $\funSubdD{F}{
    \vx} \geq \frac{\xi(\vx,\, \epsilon)}{2L_f} \geq \frac{\epsilon}{2L_tL_f}$,
    which proves the second property in \cref{def:tras:func-subd} and the last
    statement in \cref{thm:tras:func-subd:outer-comp}.

    For any $(\vx,\,\vy,\,\epsilon) \in S^2 \times \nonNegReal$, define
    $\epsilon_1 \defeq \xi(\vx,\,\epsilon)$ and $\epsilon_2 \defeq
    \gamma(\vx,\,\epsilon)$ for simplicity of notation.
    \Cref{eqn:tras:func-subd:lb} is then proven by
    \begin{align*}
        F(\vy) - F(\vx) &\geq \max_{\alpha \in \funSubdG[\epsilon_2]{h}{f(\vx)}}
            \alpha(f(\vy) - f(\vx)) - \epsilon_2 \\
        &\geq \max_{\alpha \in \funSubdG[\epsilon_2]{h}{f(\vx)}}
            \alpha\qty(
                \max_{\vg \in \funSubdG[\epsilon_1]{f}{\vx}} \T{\vg}
                (\vy - \vx) - \epsilon_1
            ) - \epsilon_2 \\
        &\geq \max_{\alpha \in \funSubdG[\epsilon_2]{h}{f(\vx)}}
            \alpha\qty(
                \max_{\vg \in \funSubdG[\epsilon_1]{f}{\vx}} \T{\vg}(\vy - \vx)
            ) - \qty(\max_{\alpha \in \funSubdG[\epsilon_2]{h}{f(\vx)}}
                \alpha\epsilon_1) - \epsilon_2 \\
        &\geq \max_{\vg_F \in \funSubdG{F}{\vx}} \T{\vg_F}(\vy - \vx) -
            \hb \epsilon_1 - \epsilon_2
        = \max_{\vg_F \in \funSubdG{F}{\vx}} \T{\vg_F}(\vy - \vx) -
            \epsilon
    \end{align*}

    Suppose $\norm{\vy - \vx}_2 \leq \funSubdD{F}{\vx}$.
    \Cref{eqn:tras:func-subd:outer-comp:xi-def} yields $\abs{f(\vy) - f(\vx)}
    \leq \funSubdD[\epsilon_2]{h}{f(\vx)}$. \Cref{eqn:tras:func-subd:ub} is thus
    proven by
    \begin{align*}
        F(\vy) - F(\vx) &\leq
            \max_{\alpha \in \funSubdG[\epsilon_2]{h}{f(\vx)}}
                \alpha(f(\vy) - f(\vx)) +
                \frac{\beta_h}{2}\abs{f(\vy) - f(\vx)}^2 \\
        &\leq \max_{\alpha \in \funSubdG[\epsilon_2]{h}{f(\vx)}}\alpha
            \qty(
                \max_{\vg \in \funSubdG[\epsilon_1]{f}{\vx}} \T{\vg}(\vy - \vx)
                + \frac{\beta_f}{2}\norm{\vy - \vx}_2^2
            ) + \frac{\beta_h L_1^2}{2}\norm{\vy - \vx}_2^2 \\
        &\leq \max_{\vg_F \in \funSubdG{F}{\vx}} \T{\vg_F}(\vy - \vx) +
            \frac{\hb\beta_f + L_1^2\beta_h}{2}\norm{\vy - \vx}_2^2 \\
        &\leq \max_{\vg_F \in \funSubdG{F}{\vx}} \T{\vg_F}(\vy - \vx) +
            \frac{\beta_F}{2}\norm{\vy - \vx}_2^2
    \end{align*}
\end{proof}
\endgroup

\begin{corollary}[Rule of outer linearity]
    \label{thm:tras:func-subd:outer-linear}
    If $F(\vx) \defeq a f(\vx) + b$ where $a \in \posReal$ and $b \in \real$,
    then $\funSubdP{F}$ is a functional subdifferential of $F(\cdot)$  with
    associated constants $(a L_f,\, a \beta_f)$ where
    \begin{align*}
        \funSubdD{F}{\vx} \defeq \funSubdD[\frac{\epsilon}{a}]{f}{\vx},\quad
        \funSubdG{F}{\vx} \defeq \condSet{a\vg}{
            \vg \in \funSubdG[\frac{\epsilon}{a}]{f}{\vx}}
    \end{align*}
\end{corollary}
\begin{proof}
    Set $h(x) \defeq ax + b$ in \cref{thm:tras:func-subd:outer-comp}.
\end{proof}

\begingroup
\newcommand{\smax}{\sigma_{\operatorname{max}}(\vA)}
\begin{proposition}[Rule of inner linearity]
    \label{thm:tras:func-subd:inner-linear}
    Let $F(\vx) \defeq f\qty(\vA \vx + \vb)$ where $f: \real^m \mapsto \real$
    is convex, $\vA \in \real^{m \times n}$, and $\vb \in \real^m$. Let $S
    \subset \real^n$ be an open convex set. Let $T \defeq \condSet{\vA \vx +
    \vb}{\vx \in S}$. Assume $\funSubdP{f}$ is a functional subdifferential of
    $f(\cdot)$ over $T$ with associated constants $(L_f,\, \beta_f)$. Then
    $\funSubdP{F}$ is a functional subdifferential of $F(\cdot)$ over $S$ with
    associated constants $\qty(\smax L_f,\, \smax^2 \beta_f)$ where
    \begin{align*}
        \funSubdD{F}{\vx} &\defeq \frac{1}{\smax} \funSubdD{f}{\vA\vx + \vb} \\
        \funSubdG{F}{\vx} &\defeq \condSet{\T{\vA}\vg}{
            \vg \in \funSubdG{f}{\vA \vx + \vb}} \\
        \smax &\defeq \sup_{\vx:\: \norm{\vx}_2 \leq 1} \norm{\vA \vx}_2
    \end{align*}
    Of note, $\smax$ is the operator norm of $\vA$ induced by the $\ell_2$-norm,
    which equals to the largest singular value of $\vA$.
\end{proposition}
\begin{proof}
    It is straightforward to verify the first two properties in
    \cref{def:tras:func-subd}.

    For any $(\vx,\,\vy,\,\epsilon) \in S^2 \times \nonNegReal$, we have
    \begin{align*}
        F(\vy) - F(\vx)
            &= f(\vA \vy + \vb) - f(\vA \vx + \vb)
            \geq \max_{\vg_0 \in \funSubdG{f}{\vA \vx + \vb}}
                \T{\vg_0}(\vA \vy - \vA \vx) - \epsilon \\
            &= \max_{\vg_0 \in \funSubdG{f}{\vA \vx + \vb}}
                \T{\qty(\T{\vA} \vg_0)}(\vy - \vx) - \epsilon
            = \max_{\vg \in \funSubdG{F}{\vx}}\T{\vg}(\vy - \vx) - \epsilon,
    \end{align*}
    which proves \cref{eqn:tras:func-subd:lb}.

    If $\norm{\vy - \vx}_2 \leq \funSubdD{F}{\vx}$, then $\norm{\vA \vy - \vA
    \vx}_2 \leq \smax \norm{\vy - \vx}_2 \leq \funSubdD{f}{\vA\vx + \vb}$. Thus
    \cref{eqn:tras:func-subd:ub} is proven by
    \begin{align*}
        F(\vy) - F(\vx)
            &= f(\vA \vy + \vb) - f(\vA \vx + \vb) \\
            &\leq  \max_{\vg_0 \in \funSubdG{f}{\vA \vx + \vb}}
                \T{\vg_0}(\vA \vy - \vA \vx) +
                \frac{\beta_f}{2}\norm{\vA \vy - \vA \vx}_2^2 \\
            &\leq \max_{\vg \in \funSubdG{F}{\vx}}\T{\vg}(\vy - \vx) +
                \frac{\smax^2 \beta_f}{2} \norm{\vy - \vx}_2^2
    \end{align*}
\end{proof}
\endgroup

\begin{proposition}[Rule of sum]
    \label{thm:tras:func-subd:sum}
    Let $F(\vx) \defeq \sum_{i=1}^m f_i(\vx)$. Then $\funSubdP{F}$ is a
    functional subdifferential of $F(\cdot)$ with associated constants $(L_F,\,
    \beta_F)$ where
    \begin{align*}
        \funSubdD{F}{\vx} &\defeq \min_{i \in \intInt{1}{m}}
            \funSubdD[\xi_i(\vx,\,\epsilon)]{f_i}{\vx},\quad
        \V{\xi}(\vx,\,\epsilon) \in \argmax_{\V{\xi} \in \epsilon \Delta_m}
            \min_{i \in \intInt{1}{m}} \funSubdD[\xi_i]{f_i}{\vx} \\
        \funSubdG{F}{\vx} &\defeq \condSet{
            \sum_{i=1}^m \vg_i}{
                \forall i \in \intInt{1}{m}:\:
                \vg_i \in \funSubdG[\xi_i(\vx,\,\epsilon)]{f_i}{\vx}} \\
        L_F &\leq m \max_{i \in \intInt{1}{m}}  L_{f_i},\quad
        \beta_F \defeq \sum_{i=1}^m \beta_{f_i}
    \end{align*}
    One can also set $\V{\xi}(\vx,\,\epsilon) = \bar{\V{\xi}} \in \epsilon
    \Delta_m$ as a constant. Then $L_F = \max_{i \in \intInt{1}{m}}
    \frac{L_{f_i}}{\bar{\xi}_i}$.
\end{proposition}
\begin{proof}
    It is straightforward to verify that $\funSubdG{F}{\vx}$ is nonempty,
    convex, closed, and bounded. By taking $\V{\xi}(\vx,\,\epsilon) =
    \frac{\epsilon}{m} \V{1}$, we have $\funSubdD{F}{\vx} \geq
    \min_{i\in\intInt{1}{m}} \funSubdD[ \frac{\epsilon}{m}]{f_i}{\vx} \geq
    \frac{\epsilon}{2} \min_{i \in \intInt{1}{m}} \frac{1}{m L_{f_i}}$. Thus the
    first two properties in \cref{def:tras:func-subd} are satisfied. One can
    also verify that
    $ \max_{\vg \in \funSubdG{F}{\vx}} \T{\vg}(\vy - \vx)
        = \sum_{i=1}^m \max_{\vg_i \in
        \funSubdG[\xi_i(\vx,\,\epsilon)]{f_i}{\vx}} \T{\vg_i}(\vy - \vx)$,
    which implies \cref{eqn:tras:func-subd:lb,eqn:tras:func-subd:ub}.
\end{proof}

\begin{corollary}[Sum of smooth and nonsmooth functions]
    \label{thm:tras:func-subd:smooth+nonsmooth}
    \\ Let $F(\vx) \defeq f(\vx) + g(\vx)$ where $f(\cdot)$ is $\beta_f$-smooth
    and $g(\cdot)$ has a functional subdifferential $\funSubdP{g}$ with
    associated constants $(L_g,\, \beta_g)$. Then $\qty(\funSubdDSym{g},\,
    \funSubdGSym{F})$ is a functional subdifferential of $F(\cdot)$ with
    associated constants $(L_g,\, \beta_f + \beta_g)$ where
    \begin{align*}
        \funSubdG{F}{\vx} &\defeq \condSet{\nabla f(\vx) + \vc}{
            \vc \in \funSubdG{g}{\vx}}
    \end{align*}
\end{corollary}
\begin{proof}
    Take \cref{thm:tras:func-subd:from-smooth} into
    \cref{thm:tras:func-subd:sum} and set $\V{\xi}(\vx,\,\epsilon) = \T{\mqty[0
    & \epsilon]}$.
\end{proof}
\begin{remark}
    \label{rmk:tras:func-subd:easy-to-compute}
    \Cref{thm:tras:func-subd:smooth+nonsmooth} is an example where the
    functional subdifferential is easier to compute than the
    relaxed/Clarke/Goldstein subdifferentials. Computing the functional
    subdifferential of $F(\vx)=f(\vx) + g(\vx)$ at $\vx_0$ needs $\nabla
    f(\vx_0)$. However, computing the relaxed subdifferential at $\vx_0$
    requires $\nabla f(\vx')$ for nonsmooth points $\vx'$ of~$g(\cdot)$, and
    computing the other two subdifferentials is more complicated.
\end{remark}

\begingroup
\newcommand{\setAct}[1][\vx]{\mathcal{A}_{\epsilon}\qty(#1)}
\newcommand{\partSG}[1][i]{\mathcal{P}_{\epsilon;\,#1}(\vx)}
\newcommand{\partSD}[1][i]{\mathcal{D}_{\epsilon;\,#1}(\vx)}
\begin{proposition}[Rule of max]
    \label{thm:tras:func-subd:max}
    Let $F(\vx) \defeq \max_{i \in \intInt{1}{m}} f_i(\vx)$. Assume
    $f_i(\cdot)$ is $L_i$-Lipschitz. Note that $L_i$ can be different from
    $L_{f_i}$, the associated constant of the functional subdifferential
    $\funSubdP{f_i}$. Then $\funSubdP{F}$ is a functional subdifferential of
    $F(\cdot)$ with associated constants $(L_F,\, \beta_F)$ where
    \begin{align*}
        \funSubdG{F}{\vx} &\defeq \setCl \setConv
            \bigcup_{i \in \setAct} \partSG \\
        \funSubdD{F}{\vx} &\defeq \min_{i \in \intInt{1}{m}}\qty(
            \begin{cases}
                \partSD + \frac{\delta_i(\vx)}{2 L_i}
                    & \text{if } i \in \setAct \\
                \frac{\delta_i(\vx)}{\max_{j \neq k} (L_j + L_k)}
                    & \text{otherwise}
            \end{cases}
        ) \\
        L_F &\leq \max_{i \in \intInt{1}{m}} \max\qty{L_i,\, L_{f_i}},\quad
        \beta_F \defeq \max_{i \in \intInt{1}{m}} \beta_{f_i} \\
        \text{where} \quad
            \delta_i(\vx) &\defeq F(\vx) - f_i(\vx),\quad
            \setAct \defeq \condSet{i}{i \in \intInt{1}{m}
                \text{ and } \delta_i(\vx) \leq \epsilon}, \\
            \partSG &\defeq \funSubdG[\epsilon - \delta_i(\vx)]{f_i}{\vx},\quad
            \partSD \defeq \min\qty{
                \funSubdD[\epsilon - \delta_i(\vx)]{f_i}{\vx},\,
                \funSubdD{f_i}{\vx}
            }
    \end{align*}
\end{proposition}
\begin{proof}
    It is straightforward to verify the first two properties in
    \cref{def:tras:func-subd}.

    One can verify that for any $(\vx,\,\vy) \in S^2$, there exist convex
    multipliers $\qty{\alpha_i}_{i \in \setAct}$ such that
    \begin{align*}
        \max_{\vg \in \funSubdG{F}{\vx}} \T{\vg}(\vy - \vx)
        &= \sum_{i \in \setAct} \alpha_i \qty(\max_{\vg_i \in \partSG}
            \T{\vg_i}(\vy - \vx)) \\
        &\leq \sum_{i \in \setAct} \alpha_i \qty(f_i(\vy) - f_i(\vx) +
            \epsilon - \delta_i(\vx)) \\
        &= \sum_{i \in \setAct} \alpha_i \qty(f_i(\vy) + \epsilon - F(\vx))
        \leq F(\vy) - F(\vx) + \epsilon,
    \end{align*}
    which proves \cref{eqn:tras:func-subd:lb}.

    Take any $\qty(\vx,\,\vy) \in S^2$ such that $\norm{\vy - \vx}_2 \leq
    \funSubdD{F}{\vx}$. To prove \cref{eqn:tras:func-subd:ub}, we consider
    the following cases for $k \in \intInt{1}{m}$. Let $\vg_k \in \argmax_{ \vg
    \in \partSG[k]} \T{\vg}(\vy - \vx)$.

    \mypara{Case 1:} $k \in \setAct$ and $\norm{\vx - \vy}_2 \leq
    \funSubdD{f_k}{\vx}$. We have
    \begin{align*}
        f_k(\vy)
        &\leq f_k(\vx) + \T{\vg_k}(\vy - \vx) +
            \frac{\beta_{f_k}}{2} \norm{\vy - \vx}_2^2 \\
        &\leq F(\vx) + \max_{\vg \in \funSubdG{F}{\vx}} \T{\vg}(\vy - \vx)
            + \frac{\beta_F}{2} \norm{\vy - \vx}_2^2
    \end{align*}

    \mypara{Case 2:} $k \in \setAct$ and $\norm{\vx - \vy}_2 >
    \funSubdD{f_k}{\vx}$. Let $\vy_0 \defeq \vx + \funSubdD{f_k}{\vx} \frac{\vy
    - \vx}{\norm{\vy - \vx}_2}$. Then $\norm{\vy_0 - \vx}_2 =
    \funSubdD{f_k}{\vx} $. Since $\norm{\vy - \vy_0}_2 = \norm{\vy - \vx}_2 -
    \funSubdD{f_k}{\vx} \leq \partSD[k] + \frac{\delta_k(\vx)}{2L_k} -
    \funSubdD{f_k}{\vx} \leq \frac{\delta_k(\vx)}{2L_k}$, we have
    \begin{align}
        \abs{f_k(\vy) - f_k(\vy_0)} \leq \frac{\delta_k(\vx)}{2},\quad
        \abs{\T{g_k}(\vy - \vy_0)} \leq \frac{\delta_k(\vx)}{2}
        \label{eqn:tras:func-subd:max:case2:1}
    \end{align}

    Let $\vh_k \in \argmax_{\vg \in \partSG[k]}
    \T{\vg}(\vy_0 - \vx)$. The property of functional subdifferential of
    $f_k(\cdot)$ implies $f_k(\vy_0) - f_k(\vx) \leq \T{\vh_k}(\vy_0 - \vx) +
    \frac{\beta_{f_k}}{2} \norm{\vy_0 - \vx}_2^2$. Since $\vy_0 - \vx$ and $\vy
    - \vx$ are collinear, we have $\vh_k = \vg_k$, yielding
    \begin{align}
        f_k(\vy_0) - f_k(\vx)
            \leq \T{\vg_k}(\vy_0 - \vx) +
                \frac{\beta_{f_k}}{2} \norm{\vy_0 - \vx}_2^2
            \leq \T{\vg_k}(\vy_0 - \vx) +
                \frac{\beta_F}{2} \norm{\vy - \vx}_2^2
        \label{eqn:tras:func-subd:max:case2:2}
    \end{align}

    Combining \cref{eqn:tras:func-subd:max:case2:1,%
    eqn:tras:func-subd:max:case2:2} yields
    \begin{align*}
        f_k(\vy) &= f_k(\vx) + \qty(f_k(\vy) - f_k(\vy_0)) +
            \qty(f_k(\vy_0) - f_k(\vx)) \\
            &\leq f_k(\vx) + \frac{\delta_k(\vx)}{2} +
                \qty(\T{\vg_k}(\vy_0 - \vy)
                    + \T{\vg_k}(\vy - \vx)
                    + \frac{\beta_F}{2} \norm{\vy - \vx}_2^2) \\
            &\leq f_k(\vx) + \frac{\delta_k(\vx)}{2} +
                \qty(\frac{\delta_k(\vx)}{2}
                    + \max_{\vg \in \funSubdG{F}{\vx}} \T{\vg}(\vy - \vx)
                    + \frac{\beta_F}{2} \norm{\vy - \vx}_2^2) \\
            &= F(\vx) + \max_{\vg \in \funSubdG{F}{\vx}} \T{\vg}(\vy - \vx)
            + \frac{\beta_F}{2} \norm{\vy - \vx}_2^2
    \end{align*}

    \mypara{Case 3:} $k \in \intInt{1}{m} \setminus \setAct$. Let $t \in
    \argmax_{i \in \intInt{1}{m}} f_i(\vx)$. Define $r_k(\vz) \defeq f_t(\vz) -
    f_k(\vz)$ for $\vz \in S$. Then $\abs{r_k(\vx) - r_k(\vy)} \leq \qty(L_t +
    L_k)\norm{\vx - \vy}_2$. We also have $r_k(\vx) = \delta_k(\vx)$. Since
    $\norm{\vy - \vx}_2 \leq \frac{\delta_k(\vx)}{L_t + L_k}$, we have
    $\abs{r_k(\vy) - r_k(\vx)} \leq \delta_k(\vx)$ and thus $r_k(\vy) \geq 0$,
    which implies $f_k(\vy) \leq f_t(\vy)$.

    Combining the above three cases proves \cref{eqn:tras:func-subd:ub}.
\end{proof}
\endgroup

\begin{proposition}[Functional subdifferential of the absolute value]
    \label{thm:tras:func-subd:abs}
    \\ Let $f(\vx) \defeq \abs{x}$ for $x \in \real$. Then $\funSubdP{f}$ is a
    functional subdifferential of $f(\cdot)$ with associated constants $(1,\, 0)$
    where
    \begin{align*}
        \funSubdD{f}{x} &\defeq
            \begin{cases}
                \abs{x} & \text{if } \abs{x} > \frac{\epsilon}{2} \\
                +\infty & \text{otherwise}
            \end{cases},\quad
        \funSubdG{f}{x} \defeq
            \begin{cases}
                \qty{\sign(x)} & \text{if } \abs{x} > \frac{\epsilon}{2} \\
                \qty[-1,\, 1] & \text{otherwise}
            \end{cases}
    \end{align*}
\end{proposition}
\begin{proof}
    It is straightforward to verify all properties in
    \cref{def:tras:func-subd}.
\end{proof}

\begin{proposition}[Functional subdifferential of the $\ell_1$-norm]
    \label{thm:tras:func-subd:l1}
    \\ Let $f(\vx) \defeq \norm{\vx}_1 = \sum_{i=1}^n \abs{x_i}$ for $\vx \in
    \real^n$. Without loss of generality, assume $\abs{x_i} \leq \abs{x_{i+1}}$
    for $i \in \intInt{1}{n-1}$. Then $\funSubdP{f}$ is a functional
    subdifferential of $f(\cdot)$ with associated constants $\qty(n,\, 0)$ where
    \begin{align*}
        \funSubdD{f}{\vx} &\defeq \begin{cases}
            \abs{x_{T+1}} & \text{if } T < n \\
            +\infty & \text{otherwise}
        \end{cases},\quad
        \funSubdG{f}{\vx} \defeq
            \condSet{\vg\in\real^n}{\vl \leq \vg \leq \vu} \\
        \text{where } &
        \begin{cases}
            l_i = -1,\,u_i = 1 & \text{if } i \leq T \\
            l_i = u_i = \sign(x_i) & \text{otherwise}
        \end{cases},\quad
        T \defeq \max \condSet{t \in \intInt{1}{n}}{
            \sum_{i=1}^t \abs{x_i} \leq
            \frac{\epsilon}{2}}
    \end{align*}
\end{proposition}
\begin{proof}
    Treat $f(\cdot)$ as the sum of absolute values of coordinates. Then apply
    \cref{thm:tras:func-subd:sum} and \cref{thm:tras:func-subd:abs}.
\end{proof}
\begin{remark}
    \label{rmk:tras:func-subd:l1}
    \Cref{thm:tras:func-subd:from-relaxed} gives another candidate of the
    functional subdifferential of the $\ell_1$ norm with better associated
    constants, but it poses a computational challenge to work with the set
    $\condSet{\vy}{\min_i \abs{y_i} = 0,\, \norm{\vy-\vx}_2 \leq
    \frac{\epsilon}{2\sqrt{n}}}$. Therefore, we use \cref{thm:tras:func-subd:l1}
    in practice.
\end{remark}

\section{The \tras{} algorithm}

\subsection{The \tras{} algorithm and its convergence analysis}
\newcommand{\deltaB}{\bar{\delta}}
\newcommand{\etaB}{\bar{\eta}}
\newcommand{\etaR}{\mathring{\eta}}
\newcommand{\etaT}{\tilde{\eta}}
\begin{algorithm}[t]
    \caption{\trasFull{} (\tras)}
    \label{algo:tras}
    \begin{algorithmic}[1]
        \footnotesize
        \Input{A convex function $f: \real^n \mapsto \real$}
        \Input{A bounded closed convex set $C \subset \real^n$}
        \Input{A functional subdifferential oracle $\funSubdG{f}{\cdot}$ over $S
            \subset \real^n$ with $C \subset S$}
        \Input{A starting point $\vx_0 \in C$}
        \Input{A sequence of functional subdifferential slack parameters
            $\epsDk[0],\, \epsDk[1],\, \ldots \in \posReal$}
        \Input{A sequence of trust region constraints $\eta_0,\, \eta_1,\,
            \ldots \in \posReal$}
        \Input{Line search parameters $\tau \in (0,\,1)$ and $\rho \in (0,\,1)$
            (default: $\tau=0.8$, $\rho=0.5$)}
        \vspace{.5em}
        \begingroup
        \everymath={\displaystyle}
        \For{$k = 0, 1, \ldots$}
            \State Define a convex set $C_k \defeq \condSet{\vd \in \real^n}{
                \vx_k + \vd \in C \text{ and } \norm{\vd}_2 \leq \eta_k}$
            \State $\vd_k \gets \argmin_{\vd \in C_k}
                \qty(\max_{\vg \in \funSubdGE{f}{\vx_k}} \T{\vg}\vd)$,
                \hspace{.5em}
                $\vg_k \gets \argmax_{\vg \in \funSubdGE{f}{\vx_k}}
                    \T{\vg}\vd_k$
                \Comment{Find a descent direction\label{algo:tras:def-dk}}
            \If{$\T{\vg_k}\vd_k \geq 0$}
                \Comment{See \cref{thm:tras:terminate}}
                \State $\vx_{k+1} \gets \vx_k$
            \Else
                \State $\lambda_k \gets 1$
                \While{$f(\vx_k + \lambda_k \vd_k) > f(\vx_k) +
                        \rho \lambda_k \T{\vg_k}\vd_k$}
                        \Comment{Backtracking line search}
                        \label{algo:tras:line-search}
                    \State $\lambda_k \gets \tau \lambda_k$
                \EndWhile
                \State $\vx_{k+1} \gets \vx_k + \lambda_k \vd_k$
            \EndIf
        \EndFor
        \endgroup
    \end{algorithmic}
\end{algorithm}

\cref{algo:tras} describes the \tras{} algorithm. Below we introduce some
notation to facilitate our analysis.
\begin{definition}
    \label{def:tras:algo-sym}
    In the analysis of \cref{algo:tras}, we use the following notation:
    \begin{itemize}
        \item $\vx^* \in \argmin_{\vx \in C} f(\vx)$ is an optimal solution.
        \item $\vz_k \defeq \vx^* - \vx_k$ is the error vector at iteration $k$.
        \item $R \defeq \sup_{(\vx,\,\vy) \in C^2} \norm{\vx - \vy}_2$ is the
            diameter of $C$.
        \item $\delta_k \defeq f(\vx_k) - f(\vx^*)$ is the optimality gap at
            iteration $k$.
        \item $\deltaB_k \defeq f(\vx_k) - f(\vx^*) - \epsDk$ is the
            optimality gap relative to $\epsDk$ at iteration $k$.
        \item $(L,\,\beta)$ are the associated constants of the functional
            subdifferential.
    \end{itemize}
\end{definition}

We first show that the line search in \cref{algo:tras} is guaranteed to
terminate.

\begin{theorem}
    \label{thm:tras:terminate}
    With the notation defined in \cref{algo:tras,def:tras:algo-sym}, at each
    iteration $k$, it holds that:
    \begin{itemize}
        \item If $\T{\vg_k}\vd_k \geq 0$, then $f(\vx_k) \leq
            f(\vx^*) + \epsDk$.
        \item If $\T{\vg_k}\vd_k < 0$, then the line search on
            \lineRef{algo:tras:line-search} terminates in $T_k$ iterations where
            \begin{align}
                T_k \leq
                \log_{\tau^{-1}}
                    \max\qty{\frac{2L\eta_k}{\epsDk},\,
                    \frac{\beta \eta_k^2}{2(1-\rho) \qty(-\T{\vg_k}\vd_k)},\,
                    1
                } + 1
                \label{eqn:tras:terminate:iter}
            \end{align}
    \end{itemize}
\end{theorem}
\begin{proof}
    Let $\vz_k \defeq \vx^* - \vx_k$. Let $\vg_k^* \in \argmax_{\vg \in
    \funSubdGE{f}{\vx_k}} \T{\vg}\vz_k$. Let $\vd_k^* \defeq \alpha_k \vz_k$
    where $\alpha_k \defeq \min\qty{\frac{\eta_k}{\norm{\vz_k}_2},\,1}$. Clearly
    $\vd_k^* \in C_k$. If $\T[*]{\vg_k}\vz_k < 0$, then
    \begin{align*}
        \T{\vg_k}\vd_k = \min_{\vd \in C_k} \max_{\vg \in \funSubdGE{f}{\vx_k}}
            \T{\vd}\vg \leq \max_{\vg \in \funSubdGE{f}{\vx_k}} \T[*]{\vd_k}\vg =
            \alpha_k \T[*]{\vg_k}\vz_k < 0
    \end{align*}
    Therefore, $\T{\vg_k}\vd_k \geq 0$ implies $\T[*]{\vg_k}\vz_k \geq 0$. With
    \cref{eqn:tras:func-subd:lb} we have $f(\vx_k) \leq f(\vx^*) -
    \T[*]{\vg_k}\vz_k + \epsDk \leq f(\vx^*) + \epsDk$ when $\T{\vg_k}\vd_k \geq
    0$, which proves the first statement in \cref{thm:tras:terminate}.

    Assuming $\T{\vg_k}\vd_k < 0$ and $\lambda_k \eta_k \leq \frac{\epsDk}{2L}
    \leq \funSubdD[\epsDk]{f}{\vx}$, a sufficient condition for the line search
    on \lineRef{algo:tras:line-search} to terminate is:
    \begingroup
    \everymath={\displaystyle}
    \renewcommand{\arraystretch}{1.8}
    \begin{align*}
    \begin{array}{lcr}
        & f(\vx_k + \lambda_k \vd_k)
            \leq f(\vx_k) + \rho \lambda_k \T{\vg_k}\vd_k & \\
        \impliedby\quad &
        \lambda_k \T{\vg_k} \vd_k + \frac{\beta}{2} \lambda_k^2 \norm{\vd_k}_2^2
            \leq \rho \lambda_k \T{\vg_k}\vd_k
            &\quad\text{Applying \cref{eqn:tras:func-subd:ub}} \\
        \impliedby &
        \frac{\beta}{2}\lambda_k^2 \eta_k^2 \leq -(1-\rho)
            \lambda_k \T{\vg_k}{\vd_k}
            &\quad\text{Applying $\norm{\vd_k}_2 \leq \eta_k$} \\
        \iff &
            \lambda_k \leq -\frac{2(1-\rho)}{\beta \eta_k^2} \T{\vg_k}\vd_k &
    \end{array}
    \end{align*}
    \endgroup
    A sufficient termination condition of \lineRef{algo:tras:line-search} is
    thus $\lambda_k \leq \min\qty{
        \frac{\epsDk}{2\eta_k L},\,
        -\frac{2(1-\rho)}{\beta \eta_k^2} \T{\vg_k}\vd_k,\,
        1
    }$. \Cref{eqn:tras:terminate:iter} then follows.
\end{proof}

Now we proceed to analyze the convergence rate of \cref{algo:tras}.
\begin{lemma}
    \label{thm:tras:lambda-eta-lb}
    With the notation defined in \cref{algo:tras,def:tras:algo-sym}, at
    iteration $k$, assume $\eta_k \leq \norm{\vz_k}_2$ and $\deltaB_k > 0$. Then
    \begin{align}
        \T{\vg_k}\vd_k \leq - \frac{\eta_k}{\norm{\vz_k}_2} \deltaB_k
        \label{eqn:tras:lambda-eta-lb:gd}
    \end{align}
    When the line search on \lineRef{algo:tras:line-search} terminates, it
    holds that:
    \begin{align}
        \lambda_k\eta_k &\geq \min\qty{\tau \etaB_k,\, \eta_k}
        \quad
        \text{where } \etaB_k \defeq \min\qty{
            \frac{\epsDk}{2L},\,
            \frac{2(1-\rho)}{\beta R} \deltaB_k
        }
    \label{eqn:tras:lambda-eta-lb}
    \end{align}
\end{lemma}
\begin{proof}
    Let $\vg_k^* \in \argmax_{\vg \in \funSubdGE{f}{\vx_k}} \min_{\vd \in C_k}
    \T{\vg}\vd$ and $\vd_k^* \in \argmin_{\vd \in C_k} \T[*]{\vg_k}\vd$. Since
    $C_k$ and $\funSubdGE{f}{\vx_k}$ are both convex and compact, we have
    $\T[*]{\vg_k} \vd_k^* = \T{\vg_k}\vd_k$ due to the minimax theorem. Since
    $\eta_k \leq \norm{\vz_k}_2$ and $C$ is convex, we have $\frac{\eta_k}{
    \norm{\vz_k}_2} \vz_k \in C_k$. We thus have $\T[*]{\vg_k} \vd_k^* \leq
    \T[*]{\vg_k}\frac{\eta_k}{\norm{\vz_k}_2} \vz_k$ due to the definition of
    $\vd_k^*$. \Cref{eqn:tras:func-subd:lb} implies $\T[*]{\vg_k} \vz_k \leq
    f(\vx^*) - f(\vx_k) + \epsDk = -\deltaB_k$. Combining these results yields
    $ \T{\vg_k}\vd_k
        = \T[*]{\vg_k}\vd^*_k
        \leq \frac{\eta_k}{ \norm{\vz_k}_2} \T[*]{\vg_k} \vz_k
        \leq -\frac{\eta_k}{\norm{\vz_k}_2} \deltaB_k
    $, which proves \cref{eqn:tras:lambda-eta-lb:gd}.

    Assuming $\lambda_k \eta_k \leq \frac{\epsDk}{2L} \leq
    \funSubdD[\epsDk]{f}{\vx}$, a sufficient condition for the line search on
    \lineRef{algo:tras:line-search} to terminate is:
    \begingroup
    \everymath={\displaystyle}
    \renewcommand{\arraystretch}{1.8}
    \begin{align*}
    \begin{array}{lcr}
        & f(\vx_k + \lambda_k \vd_k)
            \leq f(\vx_k) + \rho \lambda_k \T{\vg_k}\vd_k & \\
        \impliedby\quad &
        \lambda_k \T{\vg_k} \vd_k + \frac{\beta}{2} \lambda_k^2 \norm{\vd_k}_2^2
            \leq \rho \lambda_k \T{\vg_k}\vd_k
            &\quad\text{Applying \cref{eqn:tras:func-subd:ub}} \\
        \impliedby &
            \frac{\beta}{2}\lambda_k^2\eta_k^2
                - \qty(1-\rho) \frac{\lambda_k \eta_k}{R}\deltaB_k \leq 0  &
             \substack{
                 \text{Applying \cref{eqn:tras:lambda-eta-lb:gd}} \\
                 \text{and } \norm{\vz_k}_2 \leq R
             } \\
         \iff & \lambda_k \eta_k \leq \frac{2(1-\rho)}{\beta R} \deltaB_k &
    \end{array}
    \end{align*}
    \endgroup

    Therefore, when $\lambda_k \eta_k \leq \min\qty{
        \frac{\epsDk}{2L},\,
        \frac{2(1-\rho)}{\beta R} \deltaB_k
    } = \etaB_k$, the line search termination condition is satisfied. If $\eta_k
    < \etaB_k$, then $\lambda_k=1$ suffices. Otherwise, the search procedure can
    use at most one more iteration after $\lambda_k \leq \frac{\etaB_k}{
    \eta_k}$ is satisfied. We thus have \cref{eqn:tras:lambda-eta-lb}.
\end{proof}

\begin{lemma}
    \label{thm:tras:rate-lemma}
    With the notation defined in \cref{algo:tras,def:tras:algo-sym}, also
    assume $f(\cdot)$ is $L$-Lipschitz over $S$ (we can take $L$ to be the
    largest value of the Lipschitz constant and the associated constant of the
    functional subdifferential, which does not violate the definition of
    Lipschitz continuity or \cref{def:tras:func-subd}). Let $D \in [\delta_0,\,
    c]$ be a constant where $c$ is defined below. Given an arbitrary $\epsilon
    \in \qty(0,\,\frac{2D}{e})$ where $e$ is the base of the natural logarithm,
    set $\epsDk = \frac{\epsilon}{2}$ and $\eta_k = \frac{\epsilon}{L}$ for $k
    \in \nat$. Then
    \begin{align}
        \delta_T &\leq \frac{D}{e} + \frac{\epsilon}{2}\quad
        \text{where }
        T \defeq \left\lceil \frac{c}{\epsilon} \right\rceil,\:
        c \defeq \max\qty{
            4L,\,
            \frac{\beta R}{1-\rho}
        }\frac{R}{\rho\tau}
        \label{eqn:tras:rate-lemma}
    \end{align}
\end{lemma}
\begin{proof}
    First, note that the interval $[\delta_0,\, c]$ is non-empty since $c \geq
    LR \geq \delta_0$. We also have $\epsilon < \frac{2D}{e} < D \leq c$ and
    thus $T \geq 2$.

    If $\delta_k \leq \epsilon$ for some $k \in \intInt{0}{T}$, then $\delta_T
    \leq \delta_k \leq \epsilon \leq \frac{D}{e} + \frac{\epsilon}{2}$ since the
    objective value is non-increasing. If $\eta_k > \norm{\vz_k}_2$ for some $k
    \in \intInt{0}{T}$, then $\delta_k \leq L \norm{\vz_k}_2 \leq L \eta_k \leq
    \epsilon$, which also implies \cref{eqn:tras:rate-lemma}.

    Now let's assume $\eta_k \leq \norm{\vz_k}_2$ and $\delta_k \geq \epsilon$
    for all $k \in \intInt{0}{T}$. Then $\deltaB_k \geq \epsilon - \epsDk =
    \frac{\epsilon}{2}$. Note that $\eta_k \geq \frac{\epsDk}{2L} \geq \etaB_k$,
    which implies $\lambda_k \eta_k \geq \tau \etaB_k$ by
    \cref{thm:tras:lambda-eta-lb}. Also note that $\frac{\rho\tau}{R} \etaB_k
    \geq \frac{\rho\tau}{R} \min\qty{\frac{\epsDk}{ 2L},\,
    \frac{2(1-\rho)}{\beta R} \qty(\epsilon - \epsDk)} = \frac{\epsilon}{c}$.
    With \cref{thm:tras:lambda-eta-lb} and $\norm{\vz_k}_2 \leq R$, we have
    \begin{gather*}
        \deltaB_{k+1} - \deltaB_k
            = f(\vx_{k+1}) - f(\vx_k)
            \leq \rho \lambda_k \T{\vg_k}\vd_k
            \leq -\frac{\rho}{R} \lambda_k \eta_k \deltaB_k
            \leq -\frac{\rho\tau}{R} \etaB_k \deltaB_k
            \leq -\frac{\epsilon}{c} \deltaB_k \\
        \delta_T \leq \qty(1 - \frac{\epsilon}{c})^T \deltaB_0 + \epsDk[T]
            \leq \qty(1 - \frac{\epsilon}{c})^{\frac{c}{\epsilon}} D +
                \frac{\epsilon}{2}
            \leq \qty(\lim_{x \to 0^+} \qty(1 - x)^{\frac{1}{x}}) D
                + \frac{\epsilon}{2}
            = \frac{D}{e} + \frac{\epsilon}{2}
    \end{gather*}
\end{proof}

In \cref{thm:tras:rate-lemma}, we can set $\epsilon=\qty(1 - \frac{2}{e})D$, so
that after $T = \left\lceil \frac{ce}{(e-2)D} \right\rceil$ iterations, the
optimality gap is reduced by half to $\delta_T \leq \frac{D}{2}$. We can
recursively apply this procedure to obtain an $\bigO{\frac{1}{\epsilon}}$
algorithm, as shown in the following theorem.

\begingroup
\newcommand{\Tup}[1]{T^{\qty(#1)}}
\newcommand{\epsUp}[1]{\epsilon^{\qty(#1)}}
\newcommand{\epsTUp}[1]{\bar{\epsilon}^{\qty(#1)}}
\newcommand{\etaUp}[1]{\eta^{\qty(#1)}}
\newcommand{\TT}{T^\dag}
\newcommand{\epsT}{\epsilon^\dag}
\begin{theorem}
    \label{thm:tras:rate}
    With the notation defined in \cref{algo:tras,def:tras:algo-sym}, also
    assume $f(\cdot)$ is $L$-Lipschitz over $S$. Choose $D \in [\delta_0,\, c]$
    where $c$ is defined in \cref{thm:tras:rate-lemma}. There exist sequences
    $\qty{\eta_k}_{k \in \nat}$ and $\qty{\epsDk}_{k \in \nat}$ constructed from
    $D$ and $c$ such that in \cref{algo:tras}, for any $\epsT \in \posReal$, it
    holds that
    \begin{align*}
        \forall k\geq \TT:&\:
            f(\vx_k) - f(\vx^*) \leq \epsT \\
        \text{where } \TT &\defeq \frac{2ce}{(e-2)\epsT} +
            \left\lceil \log_2 \frac{D}{\epsT} \right\rceil
            = \bigO{\frac{1}{\epsT}}
    \end{align*}
\end{theorem}
\begin{proof}
    Let $\epsTUp{i} \defeq \frac{1}{2^i}\qty(1 - \frac{2}{e})D$ and $\Tup{i}
    \defeq \left\lceil \frac{c}{\epsTUp{i}} \right\rceil
    = \left\lceil \frac{ce}{(e-2)D}2^i \right\rceil$ for $i \in \nat$.
    For any~$i \in \nat$, define the sequences $\qty{\etaUp{i}_k}_{k \in
    \intInt{0}{\Tup{i}-1}}$ and $\qty{\epsUp{i}_k}_{k \in
    \intInt{0}{\Tup{i}-1}}$ as $\etaUp{i}_k = \frac{\epsTUp{i}}{L}$ and
    $\epsUp{i}_k = \frac{\epsTUp{i}}{2}$ as in \cref{thm:tras:rate-lemma}.
    Define the sequences $\qty{\eta_k}$ and $\qty{\epsDk}$ by concatenating the
    sequences $\qty{\etaUp{i}_k}$ and $\qty{\epsUp{i}_k}$ for $i \in \nat$.

    Let $S_0 \defeq 0$ and $S_{i+1} \defeq S_i + \Tup{i}$ for $i \in \nat$. By
    induction on $i$, it is easy to verify that $\delta_{S_i} \leq
    \frac{D}{2^i}$ using \cref{thm:tras:rate-lemma}. Let $m \defeq \left\lceil
    \log_2 \frac{D}{\epsT} \right\rceil$. We have $ \delta_{S_m} \leq \epsT$
    where
    \begin{align*}
        S_m &= \sum_{i=0}^{m-1} \Tup{i}
        \leq \sum_{i=0}^{m-1} \qty(\frac{ce}{(e-2)D} 2^i + 1)
        \leq \frac{ce}{(e-2)D} \frac{2D}{\epsT} + m
        = \TT
    \end{align*}
\end{proof}
\endgroup

\subsection{Convergence analysis for strongly convex functions}
Next we show that \tras{} converges in $\bigO{\epsilon^{-0.5}}$ iterations for
strongly convex objective functions. Recall that a function $f: \real^n \mapsto
\real$ is $\alpha$-strongly convex over $S \subset \real^n$ if and only if the
function $\vx \mapsto f(\vx) - \frac{\alpha}{2} \norm{\vx}_2^2$ is convex over
$S$, which is equivalent to $f(\vy) \geq f(\vx) + \T{\vg_x}(\vy - \vx) +
\frac{\alpha}{2} \norm{\vy - \vx}_2^2$ for $(\vx,\,\vy) \in S^2$ and $\vg_x \in
\partial f(\vx)$. We have the following property for strongly convex functions:
\begin{lemma}
    \label{thm:tras:sc:gap-bound}   
    If a function $f: \real^n \mapsto \real$ is $\alpha$-strongly convex over a
    bounded closed convex set $C \subset S$ and $\vx^* \in \argmin_{\vx \in C}
    f(\vx)$, then for any $\vx \in C$,
    \begin{align}
        \norm{\vx - \vx^*}_2 \leq \sqrt{\frac{2\qty(f(\vx) - f(\vx^*))}{\alpha}}
        \label{eqn:tras:sc:gap-bound}
    \end{align}
\end{lemma}
\begin{proof}
    By \cref{thm:subgrad-prop}, there exists $\vg \in \partial f(\vx^*)$ such
    that $\T{\vg}\qty(\vx - \vx^*) \geq 0$. \Cref{eqn:tras:sc:gap-bound} is then
    proven by the following inequality due to strong convexity:
    \begin{align*}
        f(\vx) \geq f(\vx^*) + \T{\vg}(\vx - \vx^*)
            + \frac{\alpha}{2} \norm{\vx - \vx^*}_2^2
            \geq f(\vx^*) + \frac{\alpha}{2} \norm{\vx - \vx^*}_2^2
    \end{align*}
\end{proof}

\begin{lemma}
    \label{thm:tras:sc:lambda-eta-lb}
    With the notation defined in \cref{algo:tras,def:tras:algo-sym}, also
    assume $f(\cdot)$ is $\alpha$-strongly convex and $L$-Lipschitz over $S$. At
    iteration $k$, assume $\eta_k \leq \norm{\vz_k}_2$ and $\delta_k > \epsDk$.
    Then
    \begin{align}
        \T{\vg_k}\vd_k \leq -\eta_k \sqrt{\frac{\alpha}{2}} \frac{
            \delta_k - \epsDk}{\sqrt{\delta_k}}
        \label{eqn:tras:sc:lambda-eta-lb:gd}
    \end{align}
    When the line search on \lineRef{algo:tras:line-search} terminates, it
    holds that:
    \begingroup
    \renewcommand{\arraystretch}{1.8}
    \arraycolsep=.1em
    \everymath={\displaystyle}
    \begin{align}
    \begin{array}{rl}
        \lambda_k\eta_k &\geq \min\qty{\tau \etaT_k,\, \eta_k} \\
        \text{where } \etaT_k &\defeq \min\qty{
            \frac{\epsDk}{2L},\,
            \max\qty{
                \sqrt{\frac{2\alpha}{\delta_k}},\,
                \frac{\alpha}{L}
            } \frac{(1-\rho)}{\beta}\qty(\delta_k - \epsDk)
        }
    \end{array} \label{eqn:tras:sc:lambda-eta-lb}
    \end{align}
    \endgroup
\end{lemma}
\begin{proof}
    \Cref{eqn:tras:sc:gap-bound} implies $\norm{\vz_k}_2 \leq
    \sqrt{\frac{2\delta_k}{\alpha}} \leq \sqrt{\frac{2L \norm{\vz_k}_2}{
    \alpha}}$. We thus have $\norm{\vz_k}_2 \leq \frac{2L}{\alpha}$.
    \Cref{eqn:tras:sc:lambda-eta-lb:gd} is obtained by substituting
    $\norm{\vz_k}_2 \leq \sqrt{\frac{2\delta_k}{\alpha}}$ into
    \cref{eqn:tras:lambda-eta-lb:gd}. Substituting
    $\norm{\vz_k}_2 \leq \min\qty{\frac{2L}{\alpha},\, \sqrt{\frac{
    2\delta_k}{\alpha}}}$ into the proof of \cref{thm:tras:lambda-eta-lb}
    yields \cref{eqn:tras:sc:lambda-eta-lb}.
\end{proof}

\begin{theorem}
    \label{thm:tras:sc:rate}
    With the notation defined in \cref{algo:tras,def:tras:algo-sym}, also
    assume $f(\cdot)$ is $\alpha$-strongly convex and $L$-Lipschitz over $S$.
    Let $D \in \posReal$ be any constant such that $\delta_0 \leq D$. Set
    \begin{align*}
        \eta_k = \frac{a}{(k+1)^2}&,\quad
            \epsDk = \frac{aL}{4\max\qty{k^2,\,1}} \\
        \text{where }
        a \defeq \max\qty{
            \qty(\frac{5.4}{\rho \tau c})^2 \frac{L}{2\alpha},\,
            \frac{9D}{L}
        }&,\quad
        c \defeq \min\qty{
            \frac{1}{8},\,
            \frac{5\alpha(1-\rho)}{16\beta}
            }
    \end{align*}
    Then
    \begin{align}
        \forall k \geq 1:\: \delta_k \leq \frac{aL}{k^2}
        \label{eqn:tras:sc:rate}
    \end{align}
    Equivalently, for any $\epsilon \in \posReal$, it holds that
    \begin{align*}
        \forall k \geq \sqrt{\frac{aL}{\epsilon}}:\:
            f(\vx_k) - f(\vx^*) \leq \epsilon
    \end{align*}
\end{theorem}
\begin{proof}
    Note that for $k \in \intInt{1}{3}$, $\frac{aL}{k^2} \geq D \geq \delta_0
    \geq \delta_k$ holds due to our choice of $a$.

    Now we prove by induction. Assume $k \geq 3$. If $\eta_k \geq
    \norm{\vz_k}_2$, we have $\delta_{k+1} \leq \delta_k \leq L\eta_k \leq
    \frac{aL}{(k+1)^2}$. If $\delta_k \leq \epsDk$, then $\delta_{k+1} \leq
    \delta_k \leq \frac{aL}{4k^2} = \frac{aL}{(k+1)^2} \frac{1}{4}\qty(1 +
    \frac{1}{k})^2 \leq \frac{aL}{(k+1)^2}$.

    Thus we assume $k \geq 3$, $\eta_k < \norm{\vz_k}_2$, $\delta_k > \epsDk$,
    and $\frac{aL}{(k+1)^2} < \delta_k \leq \frac{aL}{k^2}$. Substituting our
    assumptions and parameter choices into the definition of $\etaT_k$ in
    \cref{eqn:tras:sc:lambda-eta-lb} yields
    \begin{align*}
        \etaT_k &\geq \min\qty{
            \frac{a}{8k^2},\,
            \frac{\alpha(1-\rho)}{\beta L} \qty(\frac{aL}{(k+1)^2} -
            \frac{aL}{4k^2})} \\
        &= \min\qty{
            \frac{1}{8},\,
            \frac{\alpha(1-\rho)}{\beta}
                \qty(\frac{k^2}{(k+1)^2} - \frac{1}{4})
            } \frac{a}{k^2}
        \geq \min\qty{
            \frac{1}{8},\,
            \frac{5\alpha(1-\rho)}{16\beta}
            } \frac{a}{k^2}
        = \frac{ac}{k^2}
    \end{align*}
    For $k \geq 3$, we have $\frac{1}{(k+1)^2} = \frac{1}{k^2} \qty(1 -
    \frac{1}{k+1})^2 \geq \frac{9}{16k^2}$. We thus have $\eta_k =
    \frac{a}{(k+1)^2} \geq \frac{9a}{16k^2} \geq \frac{ac}{k^2} = \etaT_k$. With
    \cref{thm:tras:sc:lambda-eta-lb}, we have $\lambda_k \eta_k \geq
    \frac{ac\tau}{k^2}$ and
    \begin{align*}
        \delta_{k+1} - \delta_k \leq \rho\lambda_k\T{\vg_k}{\vd_k}
        \leq -\rho\tau\frac{ac}{k^2} \sqrt{\frac{\alpha}{2}} \frac{
            \delta_k - \epsDk}{\sqrt{\delta_k}}
    \end{align*}

    Note that $\delta_k \geq \frac{aL}{(k+1)^2} \geq \frac{9aL}{16k^2} =
    \frac{9}{4}\epsDk$ and $\delta_k \leq \frac{aL}{k^2} = 4\epsDk$, which is
    $\frac{3}{2}\sqrt{\epsDk} \leq \sqrt{\delta_k} \leq 2\sqrt{\epsDk}$.
    Consider the function $h(x) \defeq \frac{1}{x} - \frac{\epsDk}{x^3}$.
    We have $h'(x) = \frac{3\epsDk - x^2}{x^4}$, and thus the minimum of
    $f(\cdot)$ over $\qty[\frac{3}{2}\sqrt{\epsDk},\, 2\sqrt{\epsDk}]$ is
    obtained on the boundary, i.e., $h(x) \geq \frac{10}{27\sqrt{\epsDk}}
    = \frac{20}{27\sqrt{aL}}k$ for $x \in \qty[\frac{3}{2}\sqrt{\epsDk},\,
    2\sqrt{\epsDk}]$. Therefore,
    \begin{align*}
        \frac{\delta_{k+1}}{\delta_k}
        \leq 1 - \rho\tau\frac{ac}{k^2}
            \sqrt{\frac{\alpha}{2}} h\qty(\sqrt{\delta_k})
        \leq 1 - \frac{10 \rho\tau c}{27} \sqrt{\frac{2 \alpha a}{L}}\frac{1}{k}
    \end{align*}
    Consider the function $g_k(t) \defeq \frac{1 - \frac{t}{k}}
    {\frac{k^2}{(k+1)^2}} = 1 + \frac{(2-t)k^2 + (1-2t)k - t}{k^3}$. We have
    $g_k(t) < 1$ when $t \geq 2$ and $k > 0$. Therefore, $\delta_{k+1} \leq
    \frac{k^2}{(k+1)^2} g_k(t) \delta_k \leq \frac{k^2}{(k+1)^2} \delta_k \leq
    \frac{aL}{(k+1)^2}$ where $t = \frac{10 \rho\tau c}{27} \sqrt{\frac{2 \alpha
    a}{L}} \geq 2$ due to our choice of $a$.
\end{proof}

\subsection{Faster convergence with locally quadratic functional
subdifferential}
Compared to the projected gradient descent method with linear convergence for
strongly convex smooth functions, \cref{thm:tras:sc:rate} only guarantees a
$\bigO{\epsilon^{-0.5}}$ rate. The bottleneck is that
\cref{eqn:tras:sc:lambda-eta-lb} constrains $\etaT_k = \bigO{\epsDk} =
\bigO{\delta_k}$, so that \cref{eqn:tras:sc:lambda-eta-lb:gd} only decreases the
objective by $\bigO{\eta_k \sqrt{\delta_k}} = \bigO{\delta_k^{1.5}}$. We could
achieve linear convergence of \tras{} if $\etaT_k = \bigO{\sqrt{\delta_k}}$,
which intuitively means that our functional subdifferential could ``look ahead''
quadratically further. Therefore, we introduce the \emph{quadratic functional
subdifferential}, which enables faster convergence of the \tras{} algorithm.

\begin{definition}[Quadratic functional subdifferential]
    \label{def:tras:quad-func-subd}
    Let $\funSubdP{f}$ be a functional subdifferential of $f(\cdot)$ over $S$
    as defined in \cref{def:tras:func-subd}. Let $L_q \in \posReal \cup \qty{
    0^+}$ be a constant. We call $\funSubdP{f}$ a \emph{$L_q$-quadratic
    functional subdifferential} at $(\vx,\,\epsilon) \in S \times \posReal$ if $
    \funSubdD{f}{\vx} \geq \frac{\sqrt{\epsilon}}{L_q}$.
\end{definition}
\begin{remark}
If $f(\cdot)$ is smooth, then its functional subdifferential defined in
\cref{thm:tras:func-subd:from-smooth} is a $0^+$-quadratic functional
subdifferential at any $(\vx,\,\epsilon)$. For the $\ell_1$ norm,
\cref{thm:tras:func-subd:l1} gives an $L_q$-quadratic functional subdifferential
at $(\vx,\,\epsilon)$ if there exists $T$ such that $\sum_{i=1}^T \abs{x_i} \leq
\frac{\epsilon}{2}$ and $\abs{x_{T+1}} \geq \frac{\sqrt{\epsilon}}{L_q}$.
\end{remark}

\begin{lemma}
    \label{thm:tras:scqf:lambda-eta-lb}
    With the notation defined in \cref{algo:tras,def:tras:algo-sym}, also
    assume $f(\cdot)$ is $\alpha$-strongly convex over $S$. At iteration $k$,
    assume $\eta_k \leq \norm{\vz_k}_2$, $\delta_k > \epsDk$, and that
    $\funSubdP{f}$ is an $L_q$-quadratic functional subdifferential at
    $(\vx_k,\,\epsDk)$. Then \cref{eqn:tras:sc:lambda-eta-lb:gd} still holds.
    When the line search on \lineRef{algo:tras:line-search} terminates, it holds
    that:
    \begingroup
    \arraycolsep=.1em
    \renewcommand{\arraystretch}{1.8}
    \everymath={\displaystyle}
    \begin{align}
    \begin{array}{rl}
        \lambda_k\eta_k &\geq \min\qty{\tau \etaR_k,\, \eta_k} \\
        \text{where } \etaR_k &\defeq \min\qty{
            \frac{\sqrt{\epsDk}}{L_q},\,
            \frac{\sqrt{2\alpha}(1-\rho)}{\beta}
            \frac{\delta_k - \epsDk}{\sqrt{\delta_k}}
        }
    \end{array} \label{eqn:tras:scqf:lambda-eta-lb}
    \end{align}
    \endgroup
\end{lemma}
\begin{proof}
    Replace $\lambda_k \eta_k \leq \frac{\epsDk}{2L} \leq \funSubdD[\epsDk]{f}{
    \vx}$ with $\lambda_k \eta_k \leq \frac{\sqrt{\epsDk}}{L_q} \leq
    \funSubdD[\epsDk]{f}{\vx}$ and remove the Lipschitz assumption in the proof
    of \cref{thm:tras:lambda-eta-lb} and \cref{thm:tras:sc:lambda-eta-lb}.
\end{proof}

\begin{lemma}
    \label{thm:tras:ldqf:lambda-lb}   
    With the notation defined in \cref{algo:tras,def:tras:algo-sym}, also
    assume $f(\cdot)$ is $\alpha$-strongly convex over $S$. At iteration $k$,
    assume $\eta_k > \norm{\vz_k}_2$, $\delta_k > \epsDk$, and that
    $\funSubdP{f}$ is a $L_q$-quadratic functional subdifferential at
    $(\vx_k,\,\epsDk)$. Then
    \begin{align}
        \T{\vg_k}\vd_k \leq -\delta_k + \epsDk
        \label{eqn:tras:ldqf:lambda-lb:gd}
    \end{align}
    A sufficient condition for the line search on
    \lineRef{algo:tras:line-search} to terminate is:
    \begin{align}
        \lambda_k \leq \min\qty{
            \frac{\sqrt{\epsDk}}{\eta_k L_q},\,
            \frac{2(1-\rho)}{\beta \eta_k^2} \qty(\delta_k - \epsDk),\,
            1
        }
        \label{eqn:tras:ldqf:lambda-lb}
    \end{align}
\end{lemma}
\begin{proof}
    We have $\vz_k \in C_k$ since $\norm{\vz_k}_2 < \eta_k$, which implies
    \begin{align*}
        \T{\vg_k}\vd_k
            = \min_{\vd \in C_k}\max_{\vg \in \funSubdG{f}{\vx_k}} \T{\vg}\vd
            \leq \max_{\vg \in \funSubdG{f}{\vx_k}} \T{\vg}\vz_k
            \leq f(\vx^*) - f(\vx_k) + \epsDk = -\delta_k + \epsDk
    \end{align*}
    Assuming $\lambda_k \eta_k \leq \frac{\sqrt{\epsDk}}{L_q} \leq
    \funSubdD[\epsDk]{f}{\vx}$, \cref{eqn:tras:ldqf:lambda-lb} is proven
    similarly to \cref{thm:tras:terminate}.
\end{proof}

\begingroup
\newcommand{\betaT}{\tilde{\beta}}
\begin{theorem}
    \label{thm:tras:scqf:rate}  
    With the notation defined in \cref{algo:tras,def:tras:algo-sym}, also
    assume $f(\cdot)$ is $\alpha$-strongly convex over $S$. Let $D_k \in
    [\delta_k,\, +\infty)$ and $L_q \in \posReal \cup \qty{0^+}$ be two
    constants. Set
    \begin{align*}
        \epsDk &= a\gamma D_k,\quad
        \eta_k = b(1-a) \sqrt{\gamma D_k} \\
        \text{where }
        a &\defeq \qty(\frac{2b L_q}{\sqrt{4b^2L_q^2 + 1} + 1})^2,\quad
        b \defeq \frac{\sqrt{2\alpha}(1-\rho)}{\betaT},\quad
        \betaT \defeq \max\qty{\alpha,\,\beta}, \\
        \gamma &\defeq \frac{1}{1 + \rho\tau b(1-a)^2\sqrt{\frac{\alpha}{2}}}
        = \frac{1}{1 + \frac{\alpha}{\betaT} \rho(1-\rho) \tau (1-a)^2}
    \end{align*}

    If $\funSubdP{f}$ is an $L_q$-quadratic functional subdifferential at
    $(\vx_k,\,\epsDk)$, then
    \begin{align}
        \delta_{k+1} \leq \gamma D_k
        \label{eqn:tras:scqf:rate}
    \end{align}
\end{theorem}
\begin{proof}
    If $\delta_k < \gamma D_k$, then $\delta_{k+1} \leq \delta_k < \gamma D_k$,
    which yields \cref{eqn:tras:scqf:rate}. Thus we assume $ \gamma D_k <
    \delta_k \leq D_k$. Note that our choice of $a$ satisfies
    $\frac{\sqrt{a}}{L_q} = b(1-a)$ when $L_q > 0$.

    \mypara{Case 1:} $\eta_k \leq \norm{\vz_k}_2$. With
    \cref{thm:tras:scqf:lambda-eta-lb}, we have
    \begin{align*}
        \etaR_k &= \min\qty{
            \frac{\sqrt{\epsDk}}{L_q},\,
            b\frac{\delta_k - \epsDk}{\sqrt{\delta_k}}
        }
        \geq \min\qty{
            \frac{\sqrt{a\gamma D_k}}{L_q},\,
            b\frac{\gamma D_k - a \gamma D_k}{\sqrt{\gamma D_k}}
        } \\
        &= \min\qty{\frac{\sqrt{a}}{L_q},\, b(1-a)} \sqrt{\gamma D_k}
        = b(1-a) \sqrt{\gamma D_k} = \eta_k,
    \end{align*}
    which implies $\lambda_k\eta_k \geq \tau \eta_k$ with
    \cref{eqn:tras:scqf:lambda-eta-lb}. Combining this with
    \cref{eqn:tras:sc:lambda-eta-lb:gd} yields
    \begin{align*}
        \delta_{k+1} - \delta_k
            &\leq \rho \lambda_k \T{\vg_k}\vd_k
            \leq -\rho \tau b(1-a) \sqrt{\gamma D_k} \sqrt{\frac{\alpha}{2}}
            \frac{\delta_k - \epsDk}{\sqrt{\delta_k}} \\
            &\leq -\rho \tau b(1-a) \sqrt{\gamma D_k} \sqrt{\frac{\alpha}{2}}
                \frac{\gamma D_k - a \gamma D_k}{\sqrt{\gamma D_k}}
            = -\rho \tau b(1-a)^2 \sqrt{\frac{\alpha}{2}} \gamma D_k
    \end{align*}
    We thus have $\delta_{k+1} \leq \qty(1 - \rho \tau b(1-a)^2
    \sqrt{\frac{\alpha}{2}} \gamma) D_k = \gamma D_k$ due to the definition
    of~$\gamma$.

    \mypara{Case 2:} $\eta_k > \norm{\vz_k}_2$. Since
    $\frac{\sqrt{a}}{L_q} = b(1-a)$ and $\alpha \leq \betaT$, we have
    \begingroup
    \everymath={\displaystyle}
    \renewcommand{\arraystretch}{1.3}
    \arraycolsep=.1em
    \begin{gather}
        \begin{array}{c}
        \frac{\sqrt{\epsDk}}{\eta_k L_q} = \frac{\sqrt{a}}{b(1-a)L_q} = 1
        \end{array} \label{eqn:tras:scqf:rate:case2:cond1} \\
        \begin{array}{rl}
        \frac{2(1-\rho)}{\betaT \eta_k^2} \qty(\delta_k - \epsDk)
        &\geq \frac{2(1-\rho)}{\betaT b^2(1-a)^2 \gamma D_k}
            \qty(\gamma D_k - a \gamma D_k) \\
        &= \frac{2(1-\rho)}{\betaT b^2 (1-a)}
        \geq \frac{2(1-\rho)}{\betaT b^2}
        = \frac{\betaT}{\alpha (1-\rho)}
        \geq 1
        \end{array} \label{eqn:tras:scqf:rate:case2:cond2}
    \end{gather}
    \endgroup
    Substituting
    \cref{eqn:tras:scqf:rate:case2:cond1,eqn:tras:scqf:rate:case2:cond2} into
    \cref{eqn:tras:ldqf:lambda-lb} yields $\lambda_k = 1$. With
    \cref{eqn:tras:ldqf:lambda-lb:gd}, we have
    \begin{align*}
        \delta_{k+1} - \delta_k \leq \rho \lambda_k \T{\vg_k}\vd_k
        \leq -\rho(\delta_k - \epsDk)
        \leq -\rho(1 - a)\gamma D_k,
    \end{align*}
    which implies $\delta_{k+1} \leq \qty(1 - \rho(1 - a)\gamma)D_k$.
    Since $\gamma \geq \frac{1}{1 + \rho(1-a)}$, we have $\delta_{k+1} \leq
    \gamma D_k$.
\end{proof}

\begin{corollary}
    \label{thm:tras:sc-smooth:rate}
    With the notation of \cref{algo:tras,def:tras:algo-sym}, also assume
    $f(\cdot)$ is $\alpha$-strongly convex and $\beta$-smooth over $S$. Use
    \cref{thm:tras:func-subd:from-smooth} to define the functional
    subdifferential. Let $D \in \posReal$ be any constant such that $\delta_0
    \leq D$. Set
    \begin{align*}
        \epsDk = 0,&\quad
        \eta_k = b \sqrt{\gamma^k D} \\
        \text{where }
        b \defeq \frac{\sqrt{2\alpha}(1-\rho)}{\beta},&\quad
        \gamma \defeq \frac{1}{1 + \frac{\alpha}{\beta} \rho(1-\rho) \tau}
    \end{align*}

    Then
    \begin{align}
        \label{eqn:tras:sc-smooth:rate}
        \forall k \in \nat:\: \delta_k \leq \gamma^k D
    \end{align}
    Equivalently, for any $\epsilon \in \posReal$, it holds that
    \begin{align*}
        \forall k \geq \frac{\log D - \log \epsilon}{\log \gamma^{-1}}:\:
            f(\vx_k) - f(\vx^*) \leq \epsilon
    \end{align*}
\end{corollary}
\begin{proof}
    For strongly convex smooth functions, we have $\alpha \leq \beta$. With
    $L_q=0$, $a = 0$ and $\betaT = \beta$ in \cref{thm:tras:scqf:rate}, we can
    prove \cref{eqn:tras:sc-smooth:rate} by induction on $k$.
\end{proof}
\endgroup

\subsection{Almost-functional subdifferential for optimization beyond
\weakSmoothness{}}
When the objective function is only Lipschitz but not \weaklySmooth{}, a
functional subdifferential that enables efficient solutions to the minimax
problem can be harder to define. In this case, we can use the almost-functional
subdifferential, defined as the following:

\begin{definition}[Almost-functional subdifferential]
    \label{def:tras:almost-func-subd}
    With the notation defined in \cref{def:tras:func-subd}, a pair
    $\funSubdP{f}$ is called an \emph{almost-functional subdifferential} of
    $f(\cdot)$ over $S$ if all properties of \cref{def:tras:func-subd} are
    satisfied except that $\beta = \frac{\beta_0}{\epsilon}$ for a constant
    $\beta_0 \in \posReal$. We call $(L,\,\beta_0)$ the \emph{associated
    constants} of $\funSubdP{f}$.
\end{definition}

\begin{theorem}
    \label{thm:tras:almost-func-subd-rate}
    With the notation defined in \cref{algo:tras,def:tras:algo-sym}, assume
    $f(\cdot)$ is $L$-Lipschitz and $\funSubdP{f}$ is an almost-functional
    subdifferential of $f(\cdot)$ with associated constants $(L,\,\beta_0)$.
    Given $\epsilon \in \posReal$, let $T = \argmin \condSet{t \in \nat}{
    \delta_t \leq \epsilon}$. Then $T = \bigO{\epsilon^{-2}}$. If $f(\cdot)$ is
    also strongly convex, then $T = \bigO{\epsilon^{-1}}$.
\end{theorem}
\begin{proof}
    The first statement can be proven by setting $\beta =
    \frac{\beta_0}{\epsilon}$ in \cref{thm:tras:rate}.

    Now assume $f(\cdot)$ is $\alpha$-strongly convex. Similar to
    \cref{thm:tras:sc:rate}, we set $\eta_k = \frac{b}{k+1}$ and $\epsDk =
    \frac{bL}{4\max\qty{k,\,1}}$ for a constant $b$. By replacing $\beta$ with
    $\frac{ \beta_0}{\epsDk}$ in \cref{eqn:tras:sc:lambda-eta-lb}, we have
    $\etaT_k \geq \frac{\sqrt{2\alpha}(1-\rho)}{\beta_0} \frac{\delta_k -
    \epsDk}{\sqrt{\delta_k}}\epsDk$ for sufficiently large $k$, which, when
    combined with \cref{eqn:tras:sc:lambda-eta-lb:gd}, yields $\delta_{k+1} -
    \delta_k \leq -\frac{\rho\tau(1-\rho)\alpha}{\beta_0} \frac{\qty(\delta_k -
    \epsDk)^2}{\delta_k} \epsDk$. With a proper choice of $b$ and similar
    arguments to the proof of \cref{thm:tras:sc:rate}, we can prove $\delta_k
    \leq \frac{bL}{k}$, which yields $T = \bigO{\epsilon^{-1}}$.
\end{proof}

One important example of a Lipschitz but not \weaklySmooth{} function is the
$\ell_2$ norm. It has the following almost-Functional subdifferential.
\begingroup
\newcommand{\f}[1]{\norm{#1}_2}
\begin{proposition}[Almost-functional subdifferential of the $\ell_2$ norm]
    \label{thm:tras:l2:subd}
    \\ Let $f(\vx) = \f{\vx}$. For $\vx \in \real^n$ and $\epsilon \in
    \nonNegReal$, define
    \begin{align}
        \funSubdG{f}{\vx} &\defeq
        \begin{cases}
            \qty{\frac{\vx}{\f{\vx}}} &
                \text{if } \f{\vx} > \frac{\epsilon}{2} \\
            \condSet{\vg \in \real^n}{\f{\vg} \leq 1} & \text{otherwise}
        \end{cases}
        \label{eqn:tras:l2:subd}
    \end{align}
    Then $\funSubdG{f}{\vx}$ satisfies \cref{eqn:tras:func-subd:lb} and the
    following inequalities:
    \begingroup
    \everymath={\displaystyle}
    \begin{align}
        \f{\vx} > \frac{\epsilon}{2} &\implies \qty(
        \begin{array}{rl}
            \forall \vy &:\: \f{\vy - \vx} \leq \frac{\epsilon}{4}
                \implies \\
            & f(\vy) \leq f(\vx) + \max_{\vg \in \funSubdG{f}{\vx}}
                \T{\vg} (\vy - \vx) + \frac{2}{\epsilon} \norm{\vy - \vx}_2^2
        \end{array})
        \label{eqn:tras:l2:subd-ub1} \\
        \f{\vx} \leq \frac{\epsilon}{2} &\implies \qty(
            \forall \vy \in \real^n:\:
            f(\vy) \leq f(\vx) +
                \max_{\vg \in \funSubdG{f}{\vx}}\T{\vg} (\vy - \vx)
        )
        \label{eqn:tras:l2:subd-ub2}
    \end{align}
    \endgroup
    Therefore, $\funSubdP{f}$ is an almost-functional subdifferential of
    $f(\cdot)$ with associated constants $L=2$ and $\beta_0=4$ where
    $\funSubdD{f}{\vx} = \begin{cases}
        \frac{\epsilon}{4} & \text{if } \f{\vx} > \frac{\epsilon}{2} \\
        +\infty & \text{otherwise}
    \end{cases}$.
\end{proposition}
\begin{proof}
    \mypara{Case 1:} $\f{\vx} > \frac{\epsilon}{2}$. We have $\funSubdG{f}{\vx}
    = \qty{\nabla f(\vx)}$, which implies \cref{eqn:tras:func-subd:lb} due to
    convexity. Given $\vy$ such that $\f{\vy - \vx} \leq \frac{\epsilon}{4}$,
    let $\vd \defeq \vy - \vx$. Let $a \defeq \f{\vx}$, $b \defeq \f{\vd}$, and
    $c \defeq \frac{\T{\vx} \vd}{ab}$. Note that $c \in [-1,\,1]$.
    \Cref{eqn:tras:l2:subd-ub1} is then proven by
    \begin{align*}
        f(\vy) &- f(\vx) - \T{\nabla f(\vx)} (\vy - \vx)
        = \f{\vx + \vd} - \f{\vx} - \frac{\T{\vx} \vd}{\f{\vx}} \\
        &= \sqrt{a^2 + b^2 + 2abc} - a - bc
        = \frac{(1-c^2)b^2}{\sqrt{a^2 + b^2 + 2abc} + a + bc} \\
        &\leq \frac{(1-c^2)b^2}{a + bc + a + bc}
        \leq \frac{b^2}{2(a-b)}
        \leq \frac{b^2}{\epsilon - \frac{\epsilon}{2}}
        = \frac{2}{\epsilon} \norm{\vy - \vx}_2^2
    \end{align*}

    \mypara{Case 2:} $\f{\vx} \leq \frac{\epsilon}{2}$. We have $\max_{\vg \in
    \funSubdG{f}{\vx}} \T{\vg} (\vy - \vx) = \f{\vy - \vx}$.
    \Cref{eqn:tras:func-subd:lb} follows from $\f{\vy - \vx} \leq \f{\vy} +
    \f{\vx} = f(\vy) - f(\vx) + 2\f{\vx} \leq f(\vy) - f(\vx) + \epsilon$.
    \Cref{eqn:tras:l2:subd-ub2} follows from $\f{\vy - \vx} \geq \f{\vy} -
    \f{\vx} = f(\vy) - f(\vx)$.
\end{proof}
\endgroup

\subsection{An adaptive \tras{} implementation}
It is often infeasible to set the values of $\eta_k$ and $\epsDk$ based on our
previous convergence analysis since the functional subdifferential constants and
the Lipschitz constants can be challenging to compute. This section proposes
heuristic strategies to adaptively adjust the values of $\eta_k$ and $\epsDk$.

We set $\eta_k$ to be slightly larger than the maximum step length in the recent
few iterations as an estimation of the step length of the current iteration.
Formally, we set $\eta_k = \Gamma \max_{i \in\intInt{k-m_1}{k-1}} \lambda_i
\norm{\vd_i}_2$. We choose $m_1=8$ and $\Gamma = \tau^{-2}$ in our
implementation.

\begingroup
\newcommand{\epsDkn}{\epsilon_{k+1}}
Our convergence analysis sets $\epsDkn = c \delta_k = c(f(\vx_k) - f(\vx^*))$
for some constant $c \in \qty(0,\, 1)$. However, $f(\vx^*)$ is typically
unknown. Instead, we estimate $\delta_k$ from $f(\vx_{k-1}) - f(\vx_k)$.
Assuming $\delta_k = ak^{-p}$ for $a \in \posReal$ and $p \in \posReal$, then we
have
\begin{align*}
    \frac{f(\vx_{k-1}) - f(\vx_k)}{\delta_k}
        = \frac{\delta_{k-1}}{\delta_k} - 1
        = \qty(1 + \frac{1}{k-1})^p - 1 = \frac{p}{k} + o\qty(\frac{1}{k}),
\end{align*}
which suggests $\epsDkn \approx \frac{c}{p} k \qty(f(\vx_{k-1}) - f(\vx_k))$. In
order to obtain a more robust estimation, we consider the recent history and set
\begin{align*}
    \epsDkn &=  t_k \min\qty{\epsDk,\, c_k s_k},\quad
    c_{k+1} = c_k t_k,\quad
    s_k \defeq \max_{j \in \intInt{k-m_2+1}{k}} j\qty( f(\vx_{j-1}) - f(\vx_j))
\end{align*}
The scaling factor $c_k$ is initialized as $c_1 = 1$. We set $t_k \in \qty{1,\,
\mu^-,\, \mu^+}$ given parameters $\mu^- \in \qty(0,\, 1)$ and $\mu^+ > 1$. If
$\T{\vg_k}\vd_k \geq 0$, we have $\delta_k \leq \epsDk$ due to
\cref{thm:tras:terminate} and thus set $t_k = \mu^-$ to decrease $\epsDkn$. The
other case is to set $t_k = \mu^+$ if $\delta_k \geq 2\mu^+ \epsDk$. However,
since $\delta_k$ is unknown, we assume $\eta_k / \norm{\vz_k}_2$ is large enough
and consider the necessary condition $\T{\vg_k}\vd_k \leq -(\delta_k - \epsDk)
\leq -(2\mu^+ - 1)\epsDk$ (see
\cref{eqn:tras:lambda-eta-lb:gd,eqn:tras:ldqf:lambda-lb:gd}) given $\delta_k
\geq 2\mu^+ \epsDk$. If $\T{\vg_k}\vd_k \leq -(2\mu^+ - 1)\epsDk$, we further
check if using $\mu^+\epsDk$ results in a larger objective value decrease; if
so, we set $t_k = \mu^+$. Moreover, we adopt a randomization strategy to tune
$t_k$. Let $p_t \in \qty(0,\, 1)$ be a parameter and $U_k \sim U\qty(0,\, 1)$ be
a random variable. At the $k$-th iteration, if $U_k \leq p_t$, we randomly pick
$\epsDk' \in \qty{\mu^- \epsDk,\, \mu^+ \epsDk}$ and set $t_k = \epsDk'/\epsDk$
if doing so results in a larger objective value decrease. The parameter $p_t$
balances the effectiveness of the exploration against the additional computation
cost. We choose $m_2=8$, $\mu^- = 0.5$, $\mu^+=1.5$, and $p_t=0.2$ in our
implementation.
\endgroup

Another important aspect of a practical implementation is the termination
condition. We maintain a sequence $\qty{L_k}$ as the lower bounds of the
objective value:
\begin{align*}
    L_{k+1} &= \max\qty{L_k,\, f(\vx_k) - \epsDk + \Delta_k} \\
    \text{where }
    \Delta_k &\defeq \min_{\vd \in C} \max_{\vg \in \funSubdG{f}{\vx_k}}
        \T{\vg}\vd,\quad
    L_0 = -\infty
\end{align*}
One can verify $f(\vx^*) \geq L_k$. Given user-specified tolerance $\epsilon$
and maximum number of iterations $T$, the algorithm terminates when $f(\vx_k) -
L_k \leq \epsilon$ or $k \geq T$. Since $C_k \subset C$, we have $\Delta_k \leq
\T{\vg_k}\vd_k$ and thus we set $\Delta_k=-\infty$ without solving its minimax
value if $\T{\vg_k}\vd_k \leq L_k + \epsDk - f(\vx_k)$. For unconstrained
problems, we assume the optimum is within $B_{\sqrt{n}}[\vx_k]$ when $\eta_k
\leq \num{e-4}$, so we set $\Delta_k = \frac{\sqrt{n}}{\eta_k}\T{\vg_k}\vd_k$
when $\eta_k \leq \num{e-4}$ and set $\Delta_k = -\infty$ otherwise.


\subsection{Solving the \tras{} minimax problem}
Our implementation needs a user-provided callback function that returns the
solution to the minimax problem on \lineRef{algo:tras:def-dk} in
\cref{algo:tras}. We restate the problem below:
\begin{align}
    \text{The \tras{} minimax subproblem:}\quad
    \argmin_{\vd \in C_k} \qty(\max_{\vg \in \funSubdG{f}{\vx_k}} \T{\vg}\vd)
    \label{eqn:tras:minimax}
\end{align}

There are three general strategies for solving \cref{eqn:tras:minimax}:
\begin{enumerate}
    \item Directly solving the corresponding constrained optimization problem:
        \begingroup
        \everymath={\displaystyle}
        \renewcommand{\arraystretch}{1.3}
        \begin{align}
            \begin{array}{c}
                \min_{(\vd,\,u)} \: u \\
                \text{subject to }
                \vd \in C_k,\quad
                \max_{\vg \in \funSubdG{f}{\vx_k}} \T{\vg} \vd \leq u
            \end{array}
            \label{eqn:tras:minimax-direct}
        \end{align}
        \endgroup
    \item When $C$ is large enough so that $C_k = B_{\eta_k}[\V{0}]$, the
        convexity and compactness of $C_k$ and $\funSubdG{f}{\vx_k}$ imply:
        \begin{align}
            \min_{\vd \in C_k} \max_{\vg \in \funSubdG{f}{\vx_k}} \T{\vg}\vd
            = \max_{\vg \in \funSubdG{f}{\vx_k}} \min_{\vd:\: \norm{\vd}_2 \leq
                \eta_k} \T{\vg} \vd
            = -\eta_k \min_{\vg \in \funSubdG{f}{\vx_k}} \norm{\vg}_2
            \label{eqn:tras:minimax-mingrad}
        \end{align}
        Let $\vg^*$ be a solution to \cref{eqn:tras:minimax-mingrad}. Then $
        \vd^* \defeq -\eta_k \frac{\vg^*}{\norm{\vg^*}_2}$ is a solution to
        \cref{eqn:tras:minimax}.
    \item Under the same condition as above (i.e., when $C_k =
        B_{\eta_k}[\V{0}]$), we can also solve a dual form:
        \begingroup
        \everymath={\displaystyle}
        \renewcommand{\arraystretch}{1.3}
        \begin{align}
            \begin{array}{c}
                \min \norm{\vd}_2 \\
                \text{subject to } \max_{\vg \in \funSubdG{f}{\vx_k}} \T{\vg}\vd
                \leq -1
            \end{array}
            \label{eqn:tras:minimax-minstep}
        \end{align}
        \endgroup
        Let $(\vd_u^*,\,u)$ be a solution to \cref{eqn:tras:minimax-direct} and
        $\vd_v^*$ be a solution to \cref{eqn:tras:minimax-minstep}. Define $v
        \defeq \norm{\vd_v^*}_2$. We have $u=0$ if and only if
        \cref{eqn:tras:minimax-minstep} is infeasible. Otherwise,
        $\norm{\vd_u^*}_2 = \eta_k$ and $u < 0 < v$. Setting $\vd = -\frac{1}{u}
        \vd_u^*$ in \cref{eqn:tras:minimax-minstep} leads to $v \leq -\eta_k/u$.
        Similarly, setting $\vd = \frac{\eta_k}{v}\vd_v^*$ in
        \cref{eqn:tras:minimax-direct} yields $u \leq -\eta_k/v$. We thus have
        $u = v$.
\end{enumerate}

As will be shown in \cref{sec:tras:experiments}, different problems may use
different formulations for best efficiency and/or best numerical stability. Here
we consider a typical example where $C_k = B_{\eta_k}[\V{0}]$ and
$\funSubdG{f}{\vx_k}$ is a polytope with $p$ vertices. Let the columns of $\vG
\in \real^{n \times p}$ be the vertices of $\funSubdG{f}{\vx_k}$.
\Cref{eqn:tras:minimax-direct} becomes a Second-Order Cone Program (SOCP) with
the constraints $\norm{\vd}_2 \leq \eta_k$ and $\T{\vG}\vd \leq u$.
\Cref{eqn:tras:minimax-mingrad} becomes a Quadratic Program (QP) $\min_{\vx \in
\Delta_p} \norm{\vG\vx}_2$. \Cref{eqn:tras:minimax-minstep} becomes another QP
$\min_{\T{\vG}\vd \leq -1} \norm{\vd}_2$. When $p \ll n$, solving
\cref{eqn:tras:minimax-mingrad} is often more efficient than
solving~\cref{eqn:tras:minimax-minstep}. When $p$ and $n$ are comparable,
working with \cref{eqn:tras:minimax-mingrad} may be less numerically stable
compared to \cref{eqn:tras:minimax-minstep} since we need to compute $\vg =
\vG\vx$ to obtain~$\vd$.

\section{Numerical experiments}
\label{sec:tras:experiments}

\subsection{Benchmark problems}
\label{sec:tras:experiments:prob}
We perform the numerical evaluation on nine classes of nonsmooth convex problems
described below. For each problem class, we generate 50 instances. Our benchmark
set thus contains 450 test cases.

The first six problem classes are convex nonsmooth benchmark problems used in
previous work~\citep{ haarala2004new, bagirov2014introduction,
nesterov2015quasi, karmitsa2020limited}. They are defined by a single parameter
$n$, the dimension of the problem. For each of them, we generate 50 instances by
taking $n$ uniformly spaced within $\intInt{10}{n_{\max}}$. We set $n_{\max} =
5000$ for functions with sparse gradients (MAXQ, CCB3B, and SPL) and $n_{\max} =
1200$ for others (DPL, CLQ, and CCB3A).

\mypara{MAXQ}: the generalization of MAXQ in \citet{haarala2004new}.
\begin{align*}
    f(\vx) &= \max_{i\in\intInt{1}{n}}x_i^2,\quad
    x_{0i} = \begin{cases}
        i & \text{if } i \leq \frac{n}{2} \\
        -i & \text{otherwise}
    \end{cases},\quad
    \fopt = 0
\end{align*}

\mypara{DPL}: Dense Piecewise-Linear, a.k.a. the generalization of MXHILB in
\citet{ haarala2004new}. Previous nonsmooth optimization methods failed to
converge when $n \geq 1000$ \citep{ haarala2004new, bagirov2014introduction,
karmitsa2020limited}.
\begin{align*}
    f(\vx) &= \max_{i \in \intInt{1}{n}} \abs{
        \sum_{j=1}^n \frac{x_j}{i+j-1}
    },\quad
    \vx_0 = \V{1},\quad
    \fopt = 0
\end{align*}

\mypara{CLQ}: Chained LQ in \citet{haarala2004new}.
\begin{align*}
    f(\vx) &= \sum_{i=1}^{n-1} \max\qty{
        -x_i - x_{i+1},\,
        -x_i + x_{i+1} + (x_i^2 + x_{i+1}^2 - 1)
    },\\
    \vx_0 &= 0.5 \cdot \V{1},\quad
    \fopt = -(n-1)\sqrt{2}
\end{align*}

\mypara{CCB3A}: Chained CB3 I in \citet{haarala2004new}.
\begin{align*}
    f(\vx) &= \sum_{i=1}^{n-1} \max\qty{
        x_i^2 + x_{i+1}^2,\,
        (2-x_i)^2 + (2-x_{i+1})^2,\,
        2 e^{-x_i + x_{i+1}}
    }, \\
    \vx_0 &= 2 \cdot \V{1},\quad
    \fopt = 2(n-1)
\end{align*}

\mypara{CCB3B}: Chained CB3 II in \citet{haarala2004new}.
\begin{align*}
    f(\vx) &= \max\qty{
        \sum_{i=1}^{n-1} x_i^2 + x_{i+1}^2,\,
        \sum_{i=1}^{n-1} (2-x_i)^2 + (2-x_{i+1})^2,\,
        \sum_{i=1}^{n-1} 2 e^{-x_i + x_{i+1}}
    }, \\
    \vx_0 &= 2 \cdot \V{1},\quad
    \fopt = 2(n-1)
\end{align*}

\mypara{SPL}: Sparse Piecewise-Linear, the function (63) in \citet{
nesterov2015quasi}.
\begin{align*}
    f(\vx) &= \max\qty{
        \abs{x_1},\,
        \max_{i \in \intInt{2}{n}} \abs{x_i - 2 x_{i-1}}
    },\quad
    \vx_0 = \V{1},\quad
    \fopt = 0
\end{align*}

The next two problems are classic sparse linear models for regression and
classification. The parameters $m \in \posInt$, $n \in \posInt$, and $s \in
(0,\,1)$ denote the number of samples, the number of features, and the sparsity
of the solution, respectively. LLC has an additional parameter $k \in \posInt$
for the number of classes. We generate the problem instances by sampling $m$ and
$n$ uniformly from $\intInt{8}{2048}$ and $k$ uniformly from $\intInt{3}{10}$
while rejecting those with $mn > 1024^2$ for LLR and $mnk > 5\times 1024^2$ for
LLC. We sample $\lambda$ log-uniformly from $[\num{e-4},\,0.1]$. We set the
sparsity $s=0.05$. Data matrices are sampled according to the distributions
described below.

We use $\bernoulli{s}{n}$ to denote a random vector in $\real^n$ whose entries
are independent Bernoulli distributions with probability $s$. We use
$\simplexUniform{n}$ to denote the uniform distribution over $\Delta_n$, which
can be sampled by normalizing $n$ independent samples from an exponential
distribution to have a unit $\ell_1$-norm \citep{devroye2013non}.

\mypara{LLR}: Lasso Linear Regression.
\begin{align*}
    f(\vx) &= \frac{1}{2m} \norm{\vA \vx - \vb}_2^2 + \lambda \norm{\vx}_1,
    \quad \vx_0 = \V{0} \\
    \text{where  }
    & \vA \sim \normal{\V{0}_{m \times n}}{\vI},
    \quad \vb \defeq \vb_0 + \vb_N,
    \quad \vb_0 \defeq \vA (\vx_T \odot \vx_M), \\
    &\vx_T \sim \normal{\V{0}_n}{\vI},
    \quad \vx_M \sim \bernoulli{s}{n},
    \quad \vb_N \sim \normal{\V{0}_m}{\frac{0.05 \norm{\vb_0}_1}{m} \vI}
\end{align*}

\mypara{LLC}: Lasso Linear Classification.
\begin{align*}
    f(\vx) &= \frac{1}{m} \sum_{i=1}^m -\log \qty(
        \frac{\exp(\vA_i \T{\vx_{b_i}})}{\sum_{j=1}^k \exp(\vA_i \T{\vx_j})}
    ) + \lambda \norm{\vx}_1,\quad
    \vx \in \real^{k \times n},\quad
    \vx_0 = \V{0} \\
    \text{where  }&
    \vA \in \real^{m\times n},\, \vb \in \intInt{1}{k}^m
\end{align*}
To generate the data matrices for LLC, we first sample $b_i$ uniformly from
$\intInt{1}{k}$ for $i \in \intInt{1}{m}$. Then we generate a solution matrix
$\vX \in \real^{k \times n}$ by sampling each row $\vX_i$ from
$\normal{\V{0}_n}{\vI} \odot \bernoulli{s}{n}$ while rejecting if its angle with
any previous row is less than $\pi/k$. Next we generate a noisy coefficient
matrix $\vC \in \real^{m \times k}$ by setting $ C_{ib_i} = 1$ and sampling
$C_{ij}$ for $j \neq b_i$ from $0.1 \simplexUniform{k-1}$ for each $i \in
\intInt{1}{m}$. Finally, we set $\vA \defeq \bar{\vA} + \vA_N$ where $\bar{\vA}
\defeq \vC \vX$ and $\vA_N \sim \normal{\V{0}_{m \times n}}{\qty(
\frac{0.05}{mn}\sum \abs{\bar{\vA}_{ij}}) \vI}$.

The above problems are all unconstrained. Next we introduce a constrained one.

\mypara{DG}: Distance Game
\begin{gather*}
    f(\vx) = \max_{i \in \intInt{1}{m}} \sum_{j=1}^3 f^{(j)}_i(\vx),\quad
    \vx \in \Delta_n,\quad
    \vx_0 = \frac{1}{n}\V{1} \\
    \text{where  }
    f^{(1)}_i(\vx) = \abs{\T{\va_i} \vx},\quad
    f^{(2)}_i(\vx) = \norm{\vB_i \vx}_2,\quad
    f^{(3)}_i(\vx)
    \defeq \sum_k x_{k} \log \frac{x_{k} + \epsilon}{p_{ik}+\epsilon}, \\
    \va_i \in \real^n,\quad
    \vB_i \in \real^{k \times n},\quad
    \epsilon = \num{e-8},\quad
    \vp_i \in \Delta_n
\end{gather*}
The functions $f^{(j)}_i(\cdot)$ are three different convex distance functions
(note that $f^{(3)}_i(\cdot)$ is the Kullback-Leibler divergence, which is not a
metric). The parameters $m \in \posInt$ and $n \in \posInt$ control the scale of
the problem. We generate $\va_i$ and $\vB_i$ by sampling from the standard
normal distributions. We set $k = \lfloor (n+3)/4 \rfloor$ so that $\Delta_n$
contains a nonsmooth point of $f^{(2)}_i(\cdot)$ with high probability
($\prob{\min_{\vx \in \Delta_n} f^{(2)}_i(\vx) = 0}$ is $74\%$ when $n = 10$ and
$98\%$ when $n=30$ \citep{wendel1962problem}). We sample $\vp_i$ from
$\simplexUniform{n}$. Finally, we scale $\va_i$ and $\vB_i$ so that the median
values of the sets $S_j \defeq \condSet{f^{(j)}_i(\vx_0) }{i \in \intInt{1}{m}}$
are the same for $j \in \qty{1,\,2,\,3}$. We generate the problem instances by
sampling $m$ and $n$ from $\intInt{8}{1024}$ while rejecting those with $mn >
400^2$.

\subsection{Applying the \tras{} algorithm}
The functional subdifferentials of all benchmark problems can be derived in
closed form using rules given in \cref{sec:func-subd-rule}. Below we outline how
to solve the minimax problem \cref{eqn:tras:minimax} for each problem.

Three problems have closed-form solutions. The functional subdifferential of
MAXQ is a polytope for which \cref{eqn:tras:minimax-mingrad} can be written as
$-\eta_k \min_{\vp \in \Delta_d} \sqrt{\sum_{i=1}^d (a_i p_i)^2}$ where $a_i =
2x_{k_i}$ for indices $\qty{k_1,\,\cdots,\,k_d}$ defined by
\cref{thm:tras:func-subd:max}. The solution is $p_i^* = \qty(a_i^2 \sum_{j=1}^d
a_j^{-2})^{-1}$. For LLR and LLC, their functional subdifferential is a box
$\condSet{\vg}{\vu \leq \vg \leq \vv}$, which gives a solution to
\cref{eqn:tras:minimax-mingrad} as $g_i^* = \argmin_{g \in [u_i,\,v_i]}
\abs{g}$.

The functional subdifferentials of SPL, DPL, CLQ, CCB3A, and CCB3B can be
represented as polytopes or sums of polytopes. We use
\cref{eqn:tras:minimax-direct} for DPL, CLQ, and CCB3A. For SPL, we use
\cref{eqn:tras:minimax-direct} when $n > 100$ for faster speed and
\cref{eqn:tras:minimax-minstep} when $n \leq 100$ for better numerical
stability. For CCB3B, we use \cref{eqn:tras:minimax-mingrad} since the polytope
has at most three vertices.

The remaining problem is DG. We use
\cref{thm:tras:l2:subd,thm:tras:func-subd:inner-linear} to compute the
almost-functional subdifferential for $f^{(2)}_i(\cdot)$ in DG. Let $V$ be the
set of vertices of $\funSubdG{f}{\vx_k}$ for DG. We have
$
    V = \condSet{\vg_t + \T{\vM_t} \vy_t}{t \in \intInt{1}{T},\,
    \norm{\vy_t}_2 \leq 1},
$
where $T \leq 2m$, $\vM_t = \T{\vB_i}$ or $\vM_t = \V{0}$ depending on whether
$\norm{\vB_i \vx}_2 < \epsilon_t^{(2)}/2$, and $\vg_t$ is the gradient of the
smooth part of $\sum_j f^{(j)}_i(\vx)$. Note that if $\abs{\T{\va_i} \vx} <
\epsilon_t^{(1)}/2$, then two vertices are added, with $\vg_t$ containing
$\va_i$ and $\vg_{t+1}$ containing $-\va_i$. Since $\max_{\vy:\: \norm{\vy}_2
\leq 1} \T{\vd}(\vg_t + \T{\vM_t} \vy) = \T{\vg_t}\vd + \norm{\vM_t \vd}_2$
holds for any $\vd \in \real^n$, we formulate~\cref{eqn:tras:minimax-direct} as
a SOCP with constraints $\vx_k + \vd \in \Delta_n$, $\norm{\vd}_2 \leq \eta_k$,
and $\norm{\T{\vM_t} \vd}_2 \leq u - \T{\vg_t}\vd$.

\subsection{Comparison methods and implementation details}
We compare the \tras{} algorithm with the following methods:
\begin{itemize}
    \item \textbf{GD}: The projected subgradient descent as introduced in
        \cref{sec:tras:related}.
    \item \textbf{Bundle}: The proximal bundle method~\citep{
        makela2016proximal}, which is one of the fastest and most numerically
        stable methods among popular variants of bundle methods in a previous
        evaluation~\citep[Chapter 17]{ bagirov2014introduction}.
    \item \textbf{\SAtwo}: The subgradient method with double simple
        averaging, which provides convergence guarantees for the whole sequence
        of iterates (in contrast to GD that only guarantees ergodic
        convergence). \SAtwo{} has demonstrated better performance than GD on
        SPL \citep{nesterov2015quasi}. We test it on the unconstrained problems
        (problems other than DG) since applying it to the constrained case
        requires solving a nontrivial quadratic program at each iteration.
    \item \textbf{ISTA} and \textbf{FISTA}: Iterative Shrinkage-Thresholding
        Algorithm (ISTA) and Fast ISTA (FISTA) are proximal gradient methods to
        solve problems in the form $f(\vx) = g(\vx) + h(\vx)$ where $g(\cdot)$
        is smooth and $h(\cdot)$ must be simple enough to admit a closed-form
        solution to the proximal operator. FISTA achieves
        $\bigO{\epsilon^{-0.5}}$ convergence rate when $g(\cdot)$ is Lipschitz
        (possibly not strongly convex) by incorporating a momentum term. Among
        our benchmark problems, only LLR and LLC can be solved by ISTA and
        FISTA. We use the adaptive versions described in
        \citet{beck2009fast}.
\end{itemize}

Only \tras{}, Bundle, and GD are applicable to all of our benchmark problems.
For Bundle, we use the MPBNGC Fortran implementation \citep{
makela2003multiobjective} with its Julia interface~\citep{ milz2023mpbngc}. For
methods other than Bundle, we implement them in Python with numpy. \tras{}
relies on external solvers to solve the QP and SOCP problems derived from
\cref{eqn:tras:minimax}. We use the open-source Clarabel \citep{
goulart2021clarabel} solver for SOCP and the PIQP \citep{schwan2023} solver for
QP for problems other than DPL and DG. For DPL and DG, we use the commercial
Mosek solver since it is significantly faster than Clarabel on the two problems.
Of note, Clarabel generates more accurate solutions than Mosek in our
experiments.

We do not perform any problem-specific tuning for \tras{}, Bundle, ISTA, and
FISTA. We use the default hyperparameters for Bundle except that we set the
maximum line search iterations to $100$ and the maximum number of stored
subgradients to $50$. We set the step growth parameter as $1.5$ in ISTA and
FISTA. For GD and \SAtwo{}, we use the step size $\eta_k =
\frac{\eta_0}{\sqrt{k}}$ where $\eta_0=\frac{\sqrt{n}}{L}$ for SPL, DPL, and LLR
with $L$ being the Lipschitz constant, $\eta_0=0.01$ for CCB3A and CCB3B,
$\eta_0=\num{e-4}$ for DG, and $\eta_0=1$ for others.

Our benchmark environment is a Linux workstation with an AMD Ryzen Threadripper
2970WX 24-core processor and 128 GiB of RAM. We use Python 3.11.6, numpy 1.26.2,
openblas 0.3.25, Mosek 10.1.21, Clarabel 0.6.0, and PIQP 0.2.4. All methods and
external solvers use a single thread. We set the maximum number of iterations to
$50,000$ for all methods. For TRAFS and Bundle, we set the termination threshold
of solution accuracy as \num{e-6}. All methods use double-precision
floating-point numbers. The source code of experiments is available at
\url{https://github.com/jia-kai/trafs}.

\begin{table}[t]
    \centering
    \caption{Experiment results}
    \label{tab:tras:results}
    \begin{threeparttable}
        \notsotiny
        \setlength{\tabcolsep}{.5em}
        \begin{tabular}{llrrrrrrrrr}
\toprule
 \multirow[c]{2}{*}{Problem}  &  \multirow[c]{2}{*}{Method}  &  \multicolumn{3}{c}{$\epsilon=\num{e-3}$}  &  \multicolumn{3}{c}{$\epsilon=\num{e-6}$}  &  \multicolumn{3}{c}{Termination} \\
 &  & Iter\tnote{a} & Time\tnote{a} & Solved\tnote{b} & Iter\tnote{a} & Time\tnote{a} & Solved\tnote{b} & Iter\tnote{c} & Time\tnote{c} & Error\tnote{d} \\
\midrule
\multirow[c]{4}{*}{MAXQ} & \tras{} & \textbf{1.00} & \textbf{1.02} & \textbf{100\%} & \textbf{1.00} & 1.01 & \textbf{100\%} & 146 & 0.05 & \textbf{\num{1.1e-10}} \\
 & Bundle & 44.21 & 449.45 & 12\% & 5.98 & 158.94 & 4\% & 45607 & 1282.30 & \num{9.6e3} \\
 & GD & 44.39 & 7.39 & 4\% & 7.82 & \textbf{1.00} & 2\% & 50000 & 1.47 & \num{1.2e5} \\
 & \SAtwo{} & - & - & 0\% & - & - & 0\% & 50000 & 1.78 & \num{6.6e5} \\
\cmidrule(lr){1-11}
\multirow[c]{4}{*}{DPL} & \tras{} & \textbf{1.00} & \textbf{1.51} & \textbf{100\%} & \textbf{1.00} & \textbf{1.00} & \textbf{100\%} & 77 & 27.19 & \textbf{\num{7.2e-7}} \\
 & Bundle & 24.10 & 1.68 & \textbf{100\%} & - & - & 0\% & 1891 & 11.80 & \num{3.7e-5} \\
 & GD & - & - & 0\% & - & - & 0\% & 50000 & 50.02 & 0.07 \\
 & \SAtwo{} & 572.56 & 3.02 & 2\% & - & - & 0\% & 50000 & 50.24 & 0.25 \\
\cmidrule(lr){1-11}
\multirow[c]{4}{*}{CLQ} & \tras{} & \textbf{1.18} & 43.58 & \textbf{100\%} & \textbf{1.01} & \textbf{1.00} & \textbf{100\%} & 645 & 9.89 & \textbf{\num{7.6e-10}} \\
 & Bundle & 1.50 & 254.18 & \textbf{100\%} & 93.90 & 47.01 & \textbf{100\%} & 46016 & 282.72 & \num{3.7e-7} \\
 & GD & 868.94 & 45.74 & 8\% & - & - & 0\% & 50000 & 3.74 & \num{1.2e-3} \\
 & \SAtwo{} & 4.00 & \textbf{1.00} & \textbf{100\%} & - & - & 0\% & 50000 & 4.02 & \num{2.9e-5} \\
\cmidrule(lr){1-11}
\multirow[c]{4}{*}{CCB3A} & \tras{} & 1.78 & 17.21 & \textbf{100\%} & \textbf{1.00} & \textbf{1.00} & \textbf{100\%} & 137 & 2.90 & \textbf{\num{7.2e-10}} \\
 & Bundle & \textbf{1.08} & 50.63 & \textbf{100\%} & 74.29 & 30.28 & \textbf{100\%} & 47645 & 352.36 & \num{1.8e-7} \\
 & GD & 16.97 & \textbf{1.03} & \textbf{100\%} & - & - & 0\% & 50000 & 5.74 & \num{1.1e-4} \\
 & \SAtwo{} & 52.48 & 3.32 & \textbf{100\%} & - & - & 0\% & 50000 & 5.93 & \num{8.6e-5} \\
\cmidrule(lr){1-11}
\multirow[c]{4}{*}{CCB3B} & \tras{} & 1.64 & 2.52 & \textbf{100\%} & 2.10 & \textbf{1.00} & \textbf{100\%} & 95 & 0.24 & \textbf{\num{1.1e-13}} \\
 & Bundle & \textbf{1.03} & 107.50 & \textbf{100\%} & \textbf{1.00} & 19.59 & \textbf{100\%} & 39 & 1.65 & \num{3.5e-10} \\
 & GD & 15.82 & \textbf{1.04} & \textbf{100\%} & - & - & 0\% & 50000 & 6.25 & \num{6.7e-5} \\
 & \SAtwo{} & 39.49 & 2.99 & \textbf{100\%} & - & - & 0\% & 50000 & 6.63 & \num{5.5e-5} \\
\cmidrule(lr){1-11}
\multirow[c]{4}{*}{SPL} & \tras{} & \textbf{1.00} & \textbf{1.00} & \textbf{100\%} & \textbf{1.00} & \textbf{1.00} & \textbf{100\%} & 7 & 0.09 & \textbf{\num{1.3e-7}} \\
 & Bundle & 5110.32 & 3751.44 & 60\% & 3041.50 & 2058.91 & 46\% & 37536 & 1021.29 & \num{2.4e-4} \\
 & GD & - & - & 0\% & - & - & 0\% & 50000 & 14.94 & 0.21 \\
 & \SAtwo{} & - & - & 0\% & - & - & 0\% & 50000 & 15.31 & 0.06 \\
\cmidrule(lr){1-11}
\multirow[c]{6}{*}{LLR} & \tras{} & 3.50 & 13.87 & \textbf{100\%} & 3.66 & 15.10 & 90\% & 12914 & 10.94 & \num{4.5e-7}\tnote{*} \\
 & Bundle & 3.73 & 258.57 & \textbf{100\%} & 8.25 & 352.04 & 58\% & 13171 & 134.88 & \num{4.8e-6}\tnote{*} \\
 & GD & 40.67 & 15.94 & 50\% & 39.40 & 17.43 & 28\% & 50000 & 9.78 & \num{3.3e-4}\tnote{*} \\
 & \SAtwo{} & 897.55 & 327.85 & 34\% & - & - & 0\% & 50000 & 9.91 & \num{4.4e-3}\tnote{*} \\
 & ISTA & 5.62 & 4.64 & 96\% & 3.18 & 2.66 & 80\% & 46344 & 13.20 & \num{1.8e-6}\tnote{*} \\
 & FISTA & \textbf{1.04} & \textbf{1.02} & \textbf{100\%} & \textbf{1.10} & \textbf{1.04} & \textbf{100\%} & 11093 & 2.25 & \textbf{\num{2.0e-12}}\tnote{*} \\
\cmidrule(lr){1-11}
\multirow[c]{6}{*}{LLC} & \tras{} & \textbf{1.12} & 1.98 & \textbf{100\%} & \textbf{1.04} & \textbf{1.48} & \textbf{100\%} & 1412 & 7.31 & \textbf{\num{9.2e-11}}\tnote{*} \\
 & Bundle & 30.44 & 263.16 & 98\% & - & - & 0\% & 41253 & 1864.14 & \num{1.5e-4}\tnote{*} \\
 & GD & 276.10 & 126.48 & 8\% & - & - & 0\% & 50000 & 51.33 & \num{4.4e-3}\tnote{*} \\
 & \SAtwo{} & 352.78 & 227.48 & 8\% & - & - & 0\% & 50000 & 50.87 & \num{2.8e-3}\tnote{*} \\
 & ISTA & 77.68 & 37.68 & 88\% & 23.71 & 15.64 & 22\% & 50000 & 71.66 & \num{5.3e-5}\tnote{*} \\
 & FISTA & 2.66 & \textbf{1.26} & \textbf{100\%} & 3.90 & 1.53 & \textbf{100\%} & 37156 & 58.70 & \num{3.3e-10}\tnote{*} \\
\cmidrule(lr){1-11}
\multirow[c]{3}{*}{DG} & \tras{} & \textbf{1.00} & \textbf{1.00} & \textbf{100\%} & \textbf{1.00} & \textbf{1.00} & \textbf{100\%} & 180 & 16.91 & \textbf{0}\tnote{*} \\
 & Bundle & 56.83 & 9.26 & 98\% & 8.97 & 2.45 & 36\% & 34648 & 234.19 & \num{1.9e-5}\tnote{*} \\
 & GD & 435.84 & 14.20 & 60\% & - & - & 0\% & 50000 & 150.13 & \num{8.1e-4}\tnote{*} \\
\cmidrule(lr){1-11}
\multirow[c]{3}{*}{All} & \tras{} & \textbf{1.18} & \textbf{1.61} & \textbf{100.0\%} & \textbf{1.10} & \textbf{1.01} & \textbf{98.9\%} & 1735 & 8.39 & \textbf{\num{1.2e-7}}\tnote{*} \\
 & Bundle & 8.94 & 31.70 & 85.3\% & 24.25 & 39.65 & 49.3\% & 29756 & 576.15 & \num{1.1e-4}\tnote{*} \\
 & GD & 37.80 & 2.26 & 36.7\% & 30.86 & 3.00 & 3.3\% & 50000 & 32.60 & 0.01\tnote{*} \\
\bottomrule
\end{tabular}

        \begin{tablenotes}
            \item[a] Geometric mean of metrics normalized by the per-instance
                best result over instances successfully solved by the target
                method. Lower is better.
            \item[b] Proportion of problem instances successfully solved by the
                target method. Higher is better.
            \item[c] Arithmetic mean of metrics over all problem instances. Time
                is in seconds. Lower is better. Only \tras{} and Bundle support
                termination by user-defined accuracy.
            \item[d] Shifted geometric mean (see \cref{eqn:tras:exp:shm}) of the
                final solution error when the method terminates. Lower is
                better.
            \item[*] The true minimum objective value is unknown; $\fopt$ is the
                best result of the evaluated methods.
            \item[**] Bold numbers indicate the best method under the metric for
                each problem class.
        \end{tablenotes}
    \end{threeparttable}
\end{table}

\subsection{Metrics and results}
\label{sec:tras:experiments:results}
For a minimization method and a problem instance, we define its error as
\begin{align}
    E \defeq \frac{\fmeth - \fopt}{1 + \abs{\fopt}},
    \label{eqn:tras:exp:error}
\end{align}
where $\fmeth$ is the objective value achieved by the method and $\fopt$ is the
optimal objective value (listed in \cref{sec:tras:experiments:prob} for the
first six problems) or the best objective value found by any method (for LLR,
LLC, and DG).

For each problem class and each method, we evaluate the following metrics:
\begin{enumerate}
    \item For $\epsilon\in \qty{\num{e-3},\,\num{e-6}}$, we evaluate the number
        of iterations and CPU time needed to achieve $E \leq \epsilon$. We
        normalize the metrics for each problem instance by the best method on
        that instance. We then summarize the metrics of this method by computing
        the geometric mean of the normalized metrics over problem instances on
        which the method achieves $E \leq \epsilon$ . We also report the
        proportion of problem instances on which the method achieves $E \leq
        \epsilon$.
    \item We report the arithmetic mean of numbers of iterations and CPU time
        over all problem instances before the method terminates, which could be
        due to reaching the maximum number of iterations of $50,000$, satisfying
        the user-defined accuracy of $\num{e-6}$, or encountering numerical
        issues. We also report the shifted geometric mean of the final error
        $E$ over all problem instances, defined as below:
        \begin{align}
            \operatorname{SHM}(E_1,\,\cdots,\,E_n) &\defeq \exp\qty(
                \frac{1}{n}\sum_{i=1}^n \ln\qty(E_i + s)
            ) - s \label{eqn:tras:exp:shm},\quad
            \text{where } s = \num{e-6}
        \end{align}
\end{enumerate}

\Cref{tab:tras:results} presents our evaluation results. \tras{} successfully
solves all problem instances to $\epsilon=\num{e-3}$ accuracy. Under the setting
of $\epsilon=\num{e-6}$, \tras{} solves \trasSolveRate{} of the problem
instances, compared to \bundleSolveRate{} of Bundle and \gdSolveRate{} of GD. On
problems other than LLR with $\epsilon=\num{e-6}$, \tras{} is the fastest method
and successfully solves all instances; on LLR, \tras{} is slower than FISTA
which is a more specialized method with better convergence guarantees for
non-strongly convex problems. When compared to Bundle, the second-best method in
terms of the number of successfully solved instances over all problem classes,
\tras{} is \trasSpeedupA{} times faster on instances solved by both
to~$\epsilon=\num{e-3}$ and \trasSpeedupB{} times faster on instances solved by
both to~$\epsilon=\num{e-6}$.


\section{Conclusion}
This work presents the \tras{} algorithm for nonsmooth convex optimization.
\tras{} utilizes the functional subdifferential to guarantee sufficient progress
in each iteration to deliver an iteration complexity of $\bigO{\epsilon^{-1}}$
for Lipschitz functions and $\bigO{\epsilon^{-0.5}}$ for strongly convex
Lipschitz functions. These iteration complexities are better than the previously
best-known bounds of $\bigO{\epsilon^{-2}}$ and $\bigO{\epsilon^{-1}}$ in the
two settings, respectively. \tras{} assumes the ability to solve a minimax
problem involving the functional subdifferential in each iteration. We have
presented compositional rules to compute the functional subdifferential that
enable efficient solutions to the minimax problem for many functions of
practical interest. In the numerical experiments, our adaptive \tras{}
implementation achieves \trasSpeedupB{} times faster convergence and solves
twice as many problems compared to the second-best method.

\FloatBarrier

\bibliographystyle{abbrvnat}
\bibliography{refs}

\end{document}